\documentclass[12pt]{article}
\newfont{\diagram}{figfont scaled 1200}
\newcommand{\FigFunc}{1}
\newcommand{\FigRamp}{2}
\newcommand{\FigCircuit}{3}
\newcommand{\FigLoes}{4}
\begin{document}
\title{Laplace transformation updated}
\author{Ernst Terhardt (M\"unchen, Germany)}

\date{}
\maketitle
 
\begin{abstract}
The traditional theory of Laplace transformation (TLT)
as it was put forward by Gustav Doetsch was principally intended
to provide an operator calculus for ordinary
derivable functions of the $t$-domain. As TLT does not account for the
behavior of the inverse L-transform at $t=0$ its validity is
essentially confined to $t>0$. However, from solutions of linear differential
equations (DEs) one can discern
that the behavior of functions for $t\le0$ actually
is significant. It turns out that in TLT several fundamental features of
Laplace transformation (LT) are not consistently accounted for.
To get LT consistent one has to make it consistent
with the theory of Fourier transformation, and this requires that the behavior
of both the original function and of the pertinent inverse
L-transform has to be accounted for in the entire $t$-domain, i.e., for
$-\infty<t<+\infty$. When this requirement is observed there
emerges a new approach to LT which is liberated from TLT's
deficiencies and reveals certain implications of LT that
previously have either passed unnoticed or were not taken seriously.
The new approach is described; its implications are far-reaching and
heavily affect, in particular, LT's theorems for derivation/integration and
the solution of linear DEs.
\end{abstract}

\section{Introduction}
\label{Sintro}The mathematical concept which customarily is
addressed by the term {\em Laplace Transformation} (LT)
was primarily designed as a method
for the solution of linear differential equations (DEs), i.e., by a
kind of operator calculus. Though the theory of LT dates back to 
Leibniz, Euler, Laplace, Petzval, and many 20th-century authors, its presently
prevalent form, for example in \cite{Oppenheim,PhillipsParr,Weisstein},
was essentially worked out by Gustav Doetsch
\cite{Doetsch1937,Doetsch1950,Doetsch1955,Doetsch1970,Doetsch1974}.
This version of the theory is
in the present article addressed as the {\em traditional theory of LT}, TLT.

The theory of Laplace transformation is based on the pair of integral
transformations
\begin{equation}
L\{f(t)\}=F(s)=\int_0^\infty f(t){\rm e}^{-st}{\ \rm d}t;
\hbox{\ }s=\sigma+{\rm i}\omega;\hbox{\ }\sigma>0;
\label{ult}
\end{equation}
\begin{equation}
L^{-1}\{F(s)\}=\varphi(t)
={1\over2\pi{\rm i}}\int_{\sigma-{\rm i}\infty}^{\sigma+{\rm i}\infty}
F(s){\rm e}^{st}{\ \rm d}s,
\label{ilt}
\end{equation}
where $f(t)$ is presumed to be a real function of the real variable $t$.
The L-transform $F(s)$ is a complex function of the complex variable
$s=\sigma+{\rm i}\omega$. In typical applications,
$t$ denotes time and $\omega$ denotes circular frequency.

In TLT the behavior of the inverse transform is customarily characterized by 
\begin{eqnarray}
\varphi(t)&=&0\hbox{\ \ for\ }t<0,\nonumber\\
&=&f(t)\hbox{\ \ for\ }t>0;
\label{varphitlt}
\end{eqnarray}
at $t=0$ it is left undefined. Indeed,
in the realm of ordinary mathematical functions -- which originally 
was envisaged by TLT -- it is impossible to describe the inverse
L-transform's behavior at $t=0$, except if $f(0)=0$; cf.\ Sect.~\ref{Sget}. 

By contrast, the original function $f(t)$ is -- and has to be --
defined at least for $0\le t$, i.e., including $t=0$. If $f(t)$ were at
$t=0$ insufficiently defined the integral (\ref{ult}), i.e., the L-transform
$L\{f(t)\}$, would not exist.

Hence, in TLT the definition interval of the original function $f(t)$ does
not match that of the pertinent inverse transform $\varphi(t)$.
In the tradition of Doetsch's theory of LT it is widely
believed that this mismatch of definition intervals is irrelevant, provided
that application of Laplace transformation is confined to the interval $t>0$ --
which in TLT is therefore envisaged. However, this assumption is  not
tenable. For example, the solution $y(t)$ of an inhomogeneous linear
DE of the type $f[y'(t),y''(t),\ldots]=x(t)$, even when envisaged only for
$t>0$, in general depends on the behavior of $x(t)$ for $-\infty<t$
\cite{Terhardt1987}. Thus, when the solution $y(t)$ is to be obtained by LT
it is
crucial which kind of function is by $L\{x(t)\}$ actually represented in the
interval $-\infty<t$.
This fact is in Sect.~\ref{SSfail} illustrated by an example.

\parskip8pt
The most obvious symptoms of TLT's inconsistency may be listed, as follows.

\parskip8pt
\noindent
a) The mismatch of definition intervals appears to disallow
concatenation of L-trans\-forms. As the inverse transform
$\varphi(t)$ is at $t=0$ undefined it can, rigorously, not be allowed to
be subject to another L-transformation. The fact that
concatenation turns out actually to be possible does not make TLT
formally consistent in this respect. In Sect.~\ref{SSunilat} it is
explained why concatenation of L-transforms is possible.

\parskip8pt
\noindent
b) The mismatch of definition intervals appears to exclude impulse functions
at $t=0$ -- such as $\delta(t)$ -- from L-trans\-for\-ma\-tion. The
delta-impulse $\delta(t)$ is only at $t=0$ different from null; thus,
according to TLT $\delta(t)$
can not exist in the inverse L-transform. The fact that
the impulse functions actually can be retrieved from their L-transforms
does not cure this formal inconsistency of TLT.
In Sect.~\ref{SSimpulse} it is explained why impulse functions 
are LT-consistent.

\parskip8pt
\noindent
c) TLT's derivation theorem
\begin{equation}
L\{f'(t)\}=sL\{f(t)\}-f(0)
\label{dttlt}
\end{equation}
is in conflict with the definition interval of TLT, i.e., $t>0$.
The real constant $f(0)$ that in (\ref{dttlt}) appears in the L-domain
represents a $t$-domain function of its own, namely, the delta impulse
$f(0)\delta(t)$ at $t=0$. This impulse is outside TLT's definition interval
and therefore should be regarded as irrelevant \cite{Terhardt1986}. 
Inconsequently, TLT praises the ``initial value'' $f(0)$ 
as one of its most advantageous features, as that value
eventually appears in the solutions of linear DEs and is 
utilized to account for a system's initial state at $t=0$.
By its inconsequent treatment of the initial value
TLT implicitly admits that the behavior at $t=0$ of functions and
of their L-transforms actually is significant.

\parskip8pt
\noindent
d) TLT's derivation theorem (\ref{dttlt}) is inconsistent with TLT's
integration theorem
\begin{equation}
L\{f(t)\}={1\over s}L\{f'(t)\}.
\label{ittlt}
\end{equation}
The two theorems differ by the initial value $f(0)$. In TLT, this unexpected
and unexplained discrepancy is customarily tolerated. In Sect.~\ref{SSitsecond}
the relationships between the theorems for derivation and integration
are outlined and the origin of the conflict between Eqs.~(\ref{dttlt}) and
(\ref{ittlt}) is revealed.

\parskip8pt
\noindent
e) TLT does not in general keep its promise to provide the solution of
linear DEs, i.e., for $t>0$. There are discrepancies involved
which tend to be disguised by formal pseudo-consistency.
This notion just restates what was said above about solution of linear DEs.
In Sect.~\ref{SSfail} an example is described of this kind of failure.

\parskip8pt
\noindent
f) In particular, TLT's general solution of the linear DE suffers from the
so-called initial-value conflict, In TLT it has become customary to
work around this conflict by ``patching'' the original solution. The
initial-value conflict is discussed and explained in Sect.~\ref{SSthing}. 

These observations indicate that
in TLT certain fundamental aspects of LT's behavior are not consistently
accounted for.
The mismatch of definition intervals and its consequences need to be
resolved rather than 
circumvented and/or ignored. Laplace transformation needs to
become updated.

The present article offers an outline of a new, alternative 
approach to LT.
The problem of LT-consistency of $t$-domain functions is re-inspected and
resolved. From the results such obtained there emerge
new methods for the solution by LT of both the linear inhomogeneous and
the linear homogeneous DE.

Crucial results and observations are noted and emphasized as
{\em theorems}. Several
familiar theorems of LT, such as, e.g. the fundamental

\parskip8pt
\noindent
{\bf superposition theorem}
\begin{eqnarray}
L\{c_1f_1(t)+c_2f_2(t)+\ldots\}
&=&L\{c_1f_1(t)\}+L\{c_2f_2(t)\}+\ldots\nonumber\\
&=&c_1L\{f_1(t)\}+c_2L\{f_2(t)\}+\ldots,
\label{linsuper}
\end{eqnarray}
are not affected by the new insights. Another group of
theorems become restated and re-justified without assuming a new form.
A third group includes theorems that become more or less drastically
modified as compared to their familiar form. Eventually, there is a fourth
group, i.e., of theorems which may be regarded as new -- at least in so far
as in TLT they do not play a role.

Several complementary explanations and examples are exiled into an appendix.
This article is based on, and complements, earlier related work of the
present author \cite{Terhardt1985,Terhardt1986,Terhardt1987}.

\section{Getting Laplace transformation consistent}
\label{Sget}TLT's heel of Achilles lies at $t=0$.
For LT to be consistent it is not sufficient that $\varphi(t)=f(t)$ for
$t>0$; rather, the condition
\begin{equation}
\varphi(0)=f(0)
\label{varphinull}
\end{equation}
must also be fulfilled.
To suggest the implications of this requirement,
the behavior of the inverse L-transform pertinent to a continual
function $f(t)$ is illustrated in Fig.~\FigFunc.
\begin{figure}[h]
\begin{minipage}[t]{160mm}
\begin{minipage}[t]{74mm}
\hbox{\diagram\char000}
\end{minipage}
\hfill
\parbox[b]{81mm}{\small\baselineskip10pt Fig.\ \FigFunc.\
Continual function $f(t)$ (left) and the pertinent
inverse L-transform $\varphi(t)$ (right).
$\varphi(t)$ is defined for $-\infty<t$ and includes at $t=0$
the connect interval that extends from 0 to $f(0)$
\baselineskip12pt}
\end{minipage}
\end{figure}

From the figure it becomes apparent that
for functions which at $t=0$ assume a definite unique value $f(0)\ne0$
the inverse L-transform $\varphi(t)$ includes at $t=0$ an abrupt
transitional section, i.e., from 0 to $f(0)$.
As a consequence, such type of function
can not meet the criterion (\ref{varphinull}), because a {\em section}
of $\varphi(t)$ at $t=0$
can not be equal to the {\em definite unique value} $f(0)$. This kind of
non-LT-consistency applies, in particular, to the prominent class of
continual derivable functions. By contrast, if
$f(t)$ is {\em a priori} defined as a \textit{causal} function,
i.e., $f(t)=0$ for $t<0$, the criterion may actually
be met, i.e., for certain conditions which will be discussed below.

For brevity and simplicity, the discussion of LT-consistency is in the
present article focussed on the dichotomy between {\em causal functions}
and {\em d-functions}. The term {\em d-function} denotes the class of bilateral
continual functions that are in the ordinary sense derivable 
as many times as required. It should be kept in mind, in particular, that
any solution of a linear homogeneous DE of finite order is a d-function.

To achieve an explicit account of LT's behavior at $t=0$,
the definition interval must obviously be expanded from $t>0$ into $t\le0$.
As the scope of LT's formula for inverse transformation (\ref{ilt}) already
encompasses the entire $t$-domain, it is the
scope of the formula for L-transformation (\ref{ult}) which has to be
expanded. To accomplish this kind of expansion it is not required to
challenge the definitions of LT; the bilateral
scope is already implied in Eq.~(\ref{ult}).

\subsection{Laplace transformation by Fourier transformation}
\label{SSltbyft}The implications of the latter notion become apparent when one exploits
the intimate relationship that exists between Laplace transformation
and Fourier transformation.
Equation (\ref{ult}) is equivalent to the following
``unilateral'' Fourier-trans\-for\-ma\-tion formula for the 
function $f(t)\exp(-\sigma t)$, i.e.,
\begin{equation}
L\{f(t)\}=F(s)=F(\omega,\sigma)=\int_0^\infty
f(t){\rm e}^{-\sigma t}{\rm e}^{-{\rm i}\omega t}{\ \rm d}t.
\label{ultfour}
\end{equation}
However, there is not really such a thing as unilateral Fourier transformation;
Fourier transformation is inherently bilateral. The actual
{\em analysis interval} (the ``scope'') of the
transformation (\ref{ultfour}), and thus of (\ref{ult}),
is not determined by the integral's limits but by the
reciprocal of the frequency spacing ${\rm d}f={\rm d}\omega/(2\pi)$
of the corresponding Fourier-integral representation
\cite{Terhardt1985,Terhardt1987}. The actual analysis interval
of both (\ref{ult}) and (\ref{ultfour}) encompasses the entire
$t$-domain $-\infty<t<+\infty$, and the low limit $t=0$ of the integral
(\ref{ult}) corresponds to the \textit{center of} the analysis interval.
Though these implications of Fourier transformation are fairly elementary,
apprehension of them appears to be scarce; cf.\ Sect.~\ref{SSinter}.

Hence, the unilaterality of the integration interval in both
(\ref{ultfour}) and (\ref{ult}) is not equivalent to unilaterality of the
analysis interval but rather indicates
causality of the transformed function. Equation
(\ref{ultfour}) has to be equivalent to ordinary, i.e., bilateral
Fourier transformation of a causal, i.e., 
bilaterally defined function $f_{\rm c}(t)\exp(-\sigma t)$, such that
\begin{equation}
L\{f(t)\}=F(\omega,\sigma)=
\int_0^\infty f(t){\rm e}^{-\sigma t}{\rm e}^{-{\rm i}\omega t}{\ \rm d}t
=\int_{-\infty}^{+\infty}
f_{\rm c}(t){\rm e}^{-\sigma t}{\rm e}^{-{\rm i}\omega t}{\ \rm d}t.
\label{ultfourb}
\end{equation}

For LT to be consistent with Fourier transformation it is necessary
that the causal
function $f_{\rm c}(t)$ be defined in such a way that 
the two integrals in (\ref{ultfourb}) are fully equivalent.
At first sight this requirement is met, e.g., by the definition
$f_{\rm c}(t)=0$ for $t<0$; $f_{\rm c}(t)=f(t)$ for $t\ge0$.
However, this definition of the causal function
is not actually sufficient, because it leaves the
transition from $f_{\rm c}(-0)$ to $f_{\rm c}(+0)$ undefined. There
is another condition involved: The second (bilateral)
integral in (\ref{ultfourb}) has to be consistent with the pertinent
Fourier-integral representation, i.e., the ``inverse Fourier transform''
\begin{equation}
\varphi(t)=L^{-1}\{L\{f(t)\}\}
={{\rm e}^{\sigma t}\over2\pi}\int_{-\infty}^{+\infty}
F(\omega,\sigma){\rm e}^{{\rm i}\omega t}{\ \rm d}\omega
\hbox{\ \ for\ }-\infty<t.
\label{ftinv}
\end{equation}

As is suggested in (\ref{ftinv}), this expression is equivalent to (\ref{ilt}).
The Fourier-integral representation (\ref{ftinv}) of $\varphi(t)$ is 
continuously defined for $-\infty<t<+\infty$. In particular, while
$\varphi(t)$ is causal, the transition from $\varphi(-0)$ to $\varphi(+0)$
is not undefined but either includes the so-called
{\em connect interval} (see below, in particular Sect.~\ref{SSustep}),
or its derivative(s), i.e., the impulse
functions $\delta^{(n)}(t)$ ($n=0,1,\ldots$). The Fourier integral's
evident capability to represent impulse functions proves that neither the
inverse Fourier transform nor the inverse L-transform
are at $t=0$ undefined. As a consequence, for $\varphi(0)=f_{\rm c}(0)$ to hold
$f_{\rm c}(0)$ {\em must not be left insufficiently defined}.
In TLT this requirement is ignored.

Below, the consistent definition of $f_{\rm c}(t)$ and thus the consistent
expansion of LT's transformation formula
(\ref{ult}), (\ref{ultfour}) into an equivalent bilateral form, is approached
by three steps, advancing from impulse functions to the unit step
function and finally to the entire class of causal functions.

\subsection{Impulse functions}
\label{SSimpulse}The requirements just outlined for $f_{\rm c}(t)$ are {\em a priori}
met by the delta impulse and its derivatives. Thus, Eq.~(\ref{ultfourb})
applies to $f_{\rm c}(t)=\delta^{(n)}(t)$ and one obtains

\begin{equation}
L\{\delta^{(n)}(t)\}=s^n;\hbox{\ \ }n=0,1,\ldots
\label{implt}
\end{equation}

As the impulse functions are only at $t=0$ different from null,
inverse transformation is by (\ref{ftinv}) achieved for $\exp(st)=1$,
i.e., by inverse Fourier transformation of the function $({\rm i}\omega)^n$.
One obtains
\begin{equation}
L^{-1}\{L\{\delta^{(n)}(t)\}\}=\delta^{(n)}(t);\hbox{\ \ }n=0,1,\ldots
\label{impinv}
\end{equation}

It is thus established by the theory of Fourier transformation that both
the L-transform of $\delta^{(n)}(t)$ and its inverse L-transform exist,
and that the latter is identical to $\delta^{(n)}(t)$. 
The impulse functions are LT-consistent.

As $\delta^{(n)}(t)$ is defined for $-\infty<t<+\infty$, while
$\delta^{(n)}(t)=0$ for $t\ne0$,
existence of the L-transform $L\{\delta^{(n)}(t)\}$ is consistent with
the condition that for the ``unilateral'' transformations (\ref{ult}) and
(\ref{ultfour}) to apply $f(t)$ must be defined for $t\ge0$. Thus, even when
one sticks to the ``unilateral'' form of the transformation there 
is no conflict. It is only the inverse L-transform of impulse functions
that is in conflict with TLT's assumptions.

\subsection{The unit step function and the connect function}
\label{SSustep}Also the unit step function is causal and may be regarded as just another
member of the class of impulse functions. The unit step function can be
defined as the integral of the delta impulse at $t=0$, i.e., by
\begin{equation}
u(t):=\int_{-\infty}^t\delta(\tau){\ \rm d}\tau;
\hbox{\ \ }-\infty<t<+\infty.
\label{ustepintdelta}
\end{equation}

Equation (\ref{ustepintdelta}) may be regarded as an {\em implicit}
definition of
$u(t)$. When the unit step function is explicitly defined it is important
to preserve its definition at $t=0$. This can be accomplished in the form
\begin{eqnarray}
u(t)&=&0\hbox{\ \ for\ }t<0,\nonumber\\
&=&u_0(t)\hbox{\ \ for\ }t=0;\hbox{\ \ }u_0(t)=\{0\ldots1\},\nonumber\\
&=&1\hbox{\ \ for\ }t>0.
\label{ustep}
\end{eqnarray}
This definition accounts for the abrupt transition from $0$ to $1$
that occurs at $t=0$, namely, by inclusion of the so-called
{\em unit connect function},
$u_0(t)$. The pseudo-function $u_0(t)$ may be conceived of as an infinite set
$\{0\ldots1\}$ of real numbers (a distribution) that exists at $t=0$.

Inclusion of $u_0(t)$ in
$u(t)$ is indispensable not only for LT-consistency of $u(t)$ but also
for $u(t)$ to be consistent with the concept of the delta impulse.
LT-consistency of $u(t)$ requires that the condition (\ref{varphinull}) is
met, i.e., that $L^{-1}\{L\{u(t)\}\}=u(t)$ at
$t=0$. When the interval extending from
$u(-0)=0$ to $u(+0)=1$ is left undefined, fulfilment of the
condition (\ref{varphinull}) remains undecided.

The existence of the delta impulse, which is conceptualized
as the first (non-ordinary) derivative of $u(t)$, crucially depends on
existence of the connect function $u_0(t)$. With respect to (\ref{ustep})
there holds
\begin{equation}
\delta(t)=u_0'(t)=u'(t).
\label{deltaofu}
\end{equation}
Thus, indeed, inclusion of $u_0(t)$ in $u(t)$ is indispensable. 
One can not in earnest conceive the delta impulse to be the derivative of
a gap.

As the unit step function is the integral of the delta impulse, it
follows from the results depicted in Sect.~\ref{SSimpulse} that
the unit step function as defined by
(\ref{ustepintdelta}) and (\ref{ustep}) is indeed LT-consistent;
i.e., there holds
\begin{equation}
L^{-1}\{L\{u(t)\}\}=u(t)\hbox{\ \ for\ }-\infty<t<+\infty.
\label{ustepinv}
\end{equation}

Notably, the L-transform of the unit step function is equivocal. There holds
\begin{equation}
L\{u(t)\}=L\{1\}=1/s.
\label{ustepident}
\end{equation}
With respect to the interval $t\ge0$ the functions
$f(t)=u(t)$ and $f(t)=1$ differ by the connect function $u_0(t)$ which
is included in $f(t)=u(t)$ but not in $f(t)=1$. The connect function does
not become explicitly apparent in the L-transform. As a consequence,
from the L-transform $1/s$ one can not tell whether it was obtained from
$f(t)=u(t)$ or from $f(t)=1$.

From the fact that both $u(t)$ and its derivatives (the impulse functions)
are LT-consistent one can conclude that $u_0(t)$ invariably is contained in
the \textit{inverse} L-transform. According to
(\ref{ustepinv}, \ref{ustepident}) there holds
\begin{equation}
L^{-1}\{L\{u(t)\}\}=L^{-1}\{L\{1\}\}=u(t),
\label{ustepidentinv}
\end{equation}
and $u(t)$ contains $u_0(t)$.
Thus, for a real function $f(t)$ to be LT-consistent
it is indispensable that it contains the connect function. 
This is why $f(t)=u(t)$ is LT-consistent whereas $f(t)=1$ is not.

The somewhat confusing behavior of $u_0(t)$ emerges from the fact that
$u_0(t)$ is a {\em null-function}, i.e.,
\begin{equation}
\int_{-T}^{+T}u_0(t){\ \rm d}t=0;\hbox{\ \ }T\ge0;
\label{unullint}
\end{equation}
and, therefore,
\begin{equation}
L\{u_0(t)\}=0.
\label{unulllt}
\end{equation}

Whereas $u_0(t)$ as such does not become apparent in L-transforms,
its derivatives do. The derivatives
\begin{equation}
u_0^{(n)}(t)=u^{(n)}(t)=\delta^{(n-1)}(t);\hbox{\ \ }n=1,2,\ldots
\label{unullderiv}
\end{equation}
have the L-transforms
\begin{equation}
L\{u_0^{(n)}(t)\}=L\{u^{(n)}(t)\}=L\{\delta^{(n-1)}(t)\}=s^{n-1},
\hbox{\ \ }n=1,2,\ldots
\label{deltanlt}
\end{equation}
and these are different from null.

The occurrence of the connect function in the inverse L-transform
was implicitly pointed out, e.g., by Doetsch
\cite{Doetsch1950,Doetsch1970,Doetsch1974}. He
proved that the inverse L-transform is identical to the original function
{\em except for a null-function}, i.e., a function whose integral is null.
Thus, Doetsch in effect anticipated the involvement of the connect function.
However, the fact that the null function does not become apparent in
L-transforms led him to define inverse L-transforms that differ only by
a null-function to be identical. As a consequence,
in TLT the connect function is being ignored. This is a serious mistake,
i.e., for the following reasons.

\parskip8pt\noindent
a) Without inclusion of $u_0(t)$ one can not obtain a solid
\textit{conceptual} definition of LT-consistency, as was pointed out above;
cf.\ Fig.~\FigFunc.

\parskip8pt\noindent
b) Without inclusion of $u_0(t)$ one can not obtain a solid \textit{formal}
definition of LT-con\-sis\-ten\-cy, because the \textit{unit-step
re\-dun\-dan\-cy theorem}
$u(t)\cdot u^{(n)}(t)=u^{(n)}(t)$
does not hold; cf.\ Sects.~\ref{SSunilat}, \ref{SSusred}.

\parskip8pt\noindent
c) Without inclusion of $u_0(t)$ neither the derivatives of $u(t)$, i.e.,
the impulse functions $\delta^{(n)}(t)$, nor their L-transforms
are sufficiently defined.

\subsection{Causal functions}
\label{SSunilat}Utilizing the consistent definition (\ref{ustepintdelta}, \ref{ustep})
of the unit step function, finally the causal type of function
$f_{\rm c}(t)$ can be defined in the familiar way, i.e.,
\begin{equation}
f_{\rm c}(t)=u(t)f(t).
\label{fu}
\end{equation}

It is inclusion of the unit connect function $u_0(t)$ in the unit step
function $u(t)$ that makes the expression
(\ref{fu}) \textit{universal}. Equation (\ref{fu}) holds irrespective of
whether or not
$f(t)$ itself is causal. $f(t)$ may be either a d-function $f_{\rm d}(t)$;
a causal function of the form $[u(t)f_{\rm d}(t)]$; or an impulse
function $\delta^{(n)}(t)$. This follows from the identities outlined in
Sect.~\ref{SSusred} which can be subsumed by the

\parskip12pt\noindent
{\bf unit-step redundancy theorem}
\begin{equation}
u(t)\cdot u^{(n)}(t)=u^{(n)}(t);\hbox{\ \ }n=0,1,\ldots
\label{ustepequiv}
\end{equation}
By this theorem Eq.~(\ref{fu}) remains in effect unchanged when
both sides are multiplied by $u(t)$; in particular, there holds
\begin{equation}
u(t)f_{\rm c}(t)=f_{\rm c}(t).
\label{fcred}
\end{equation}
\textit{Multiplication by $u(t)$ of any kind of causal function is redundant.}

On the basis of these observations
the bilateral formula for L-transformation envisaged in (\ref{ultfourb})
can be restated, and its equivalence to (\ref{ult}) can be expressed by the

\parskip8pt\noindent
{\bf bilaterality theorem}
\begin{equation}
L\{f(t)\}=\int_{-\infty}^{+\infty}u(t)f(t){\rm e}^{-st}{\ \rm d}t.
\label{bult}
\end{equation}
It is the LT-consistent definition (\ref{ustep}) of $u(t)$ which renders
Eq.~(\ref{bult}) consistent, i.e., by consistency with Fourier transformation;
and Eq.~(\ref{ult}) becomes consistent with Fourier transformation
by its equivalence to (\ref{bult}). The mathematical 
implications of (\ref{ult}) are the same as those of (\ref{bult}).

While the above three-step approach to the formula (\ref{bult}) is helpful
by its elucidating implications,
it should be noticed that formally the expansion of (\ref{ult}) into
(\ref{bult}) can be obtained by one single step, namely,
\begin{equation}
L\{f(t)\}=\int_0^\infty f(t){\rm e}^{-st}{\ \rm d}t
=\lim_{T\to0}\int_{-\infty}^{+\infty}r(t,T)f(t){\rm e}^{-st}{\ \rm d}t,
\label{bultlim}
\end{equation}
where $r(t,T)$ denotes the unit ramp function of which an example is
illustrated in Fig.~\FigRamp.
\begin{figure}[h]
\begin{minipage}[t]{160mm}
\begin{minipage}[t]{37mm}
\hbox{\diagram\char001}
\end{minipage}
\hfill
\parbox[b]{119mm}{\small\baselineskip10pt Fig.\ \FigRamp.\
The unit ramp function $r(t)$ as an approximation to the unit step function
$u(t)$. For $T\to0$ $u(t)$ emerges from $r(t)$, and the unit connect
function $u_0(t)$ emerges from $r_0(t)$
\baselineskip12pt}
\end{minipage}
\end{figure}

By letting $T\to0$ the unit step function $u(t)$ emerges from the unit ramp
function $r(t,T)$, and the connect function $u_0(t)$ 
emerges quite naturally from the ascending part $r_0(t)$ of the unit ramp
function.

From this approach it becomes apparent that
the transition from $u(-0)f(t)=0$ to $u(+0)f(t)=f(+0)$ is not empty
but is a vertically ascending continuous function, i.e., $u_0(t)$.
From this notion there emerges another
\textit{definition} of $u_0(t)$, namely,
\begin{equation}
u_0(t)=\lim_{T\to0}r_0(t,T).
\label{unulldef}
\end{equation}
The function $r_0(t)$ does not necessarily have to be linear -- such as
in Fig.~\FigRamp\ -- but may be represented by any kind of real continuous
function that in the interval $0\le T$ rises monotonically from 0 to 1.

\parskip8pt
Comparison of (\ref{ult}) to (\ref{bult}) reveals a crucial
feature of LT which in TLT is ignored, namely, LT's inherent
ambivalence. There holds the

\parskip8pt\noindent
{\bf ambivalence theorem}
\begin{equation}
L\{f(t)\}=L\{u(t)f(t)\}=F(s).
\label{identity}
\end{equation}
From an L-transform $F(s)$ one can not tell 
whether it was obtained from $f(t)$ or from $u(t)f(t)$.
This kind of ambivalence corresponds to the ambivalence of $L\{u(t)\}$
that was noted by (\ref{ustepident}). In both cases the ambivalence is
significant by the presence of $u_0(t)$ in $u(t)$ and by the fact that
$L\{u_0^{(n)}(t)\}=s^{n-1}\ne0$ ($n=1,2,\ldots$).

As a consequence of the equivalence of (\ref{bult}) to the formula for
Fourier-trans\-for\-ma\-tion, the inverse L-transform obtained by (\ref{ilt})
is invariably and unequivocally identical to $u(t)f(t)$. There holds the

\parskip8pt\noindent
{\bf causality theorem}  
\begin{equation}
\varphi(t)=L^{-1}\{L\{f(t)\}\}=L^{-1}\{L\{u(t)f(t)\}\}=u(t)f(t)
\hbox{\ \ for\ }-\infty<t<+\infty.
\label{varphi}
\end{equation}
The inverse L-transform is the {\em causal companion} of
the original function $f(t)$.

The latter two theorems warrant con\-sis\-ten\-cy of
\textit{con\-ca\-te\-na\-tion}
of L-trans\-for\-ma\-tions. From (\ref{identity}) and (\ref{varphi}) there
follows the

\parskip8pt\noindent
{\bf concatenation theorem}  
\begin{equation}
L\{\varphi(t)\}=L\{f(t)\}.
\label{concat}
\end{equation}
In TLT this theorem can not exist because $\varphi(0)$ is undefined
such that, rigorously, $L\{\varphi(t)\}$ can not be supposed to exist.

As the behavior of the L-transform L\{$f(t)\}$ is in every respect
characterized by the causal inverse L-transform $\varphi(t)$,
in a sense the  L-transform of $f(t)$ \textit{actually} constitutes
the L-transform of the causal function $u(t)f(t)$. 
In many contexts -- e.g., that of the
derivation/integration theorem (cf.\ Sect.\ \ref{Stheo}) --
it is indeed quite helpful to observe the

\parskip12pt\noindent
{\bf alias theorem:}

The notation $L\{f(t)\}$ of an L-transform ultimately is
an alias for $L\{u(t)f(t)\}$. 

From the above insights there emerges a concise proof of the 

\parskip8pt\noindent
{\bf convolution theorem}
\begin{equation}
L\{f_1(t)\}\cdot L\{f_2(t)\}
=L\Bigl\lbrace\int_0^tf_1(t-\tau)f_2(\tau){\ \rm d}\tau\Bigr\rbrace.
\label{convtheo}
\end{equation}
The proof of (\ref{convtheo}) can be based on the identity
\begin{equation} 
\int_{-\infty}^Tg_1(t){\ \rm d}t\cdot\int_{-\infty}^Tg_2(t){\ \rm d}t
=u(t)\int_0^T\int_0^tg_1(t-\tau)g_2(\tau){\ \rm d}\tau{\rm\ d}t;
\hbox{\ \ }T>0;\hbox{\ }t\le T,
\label{gprod}
\end{equation}
which holds for any pair of \textit{causal} functions $g_{1,2}(t)$.
Letting $g_{1,2}(t)=u(t)f_{1,2}(t)\exp(-st)$ and $T\to\infty$,
one obtains from (\ref{gprod})
\begin{equation}
\int_0^\infty f_1(t){\rm e}^{-st}{\ \rm d}t
\cdot\int_0^\infty f_2(t){\rm e}^{-st}{\ \rm d}t
=\int_0^\infty{\rm e}^{-st}\int_0^t f_1(t-\tau)f_2(\tau)
{\ \rm d}\tau{\rm\ d}t,
\label{fprod}
\end{equation}
where the ``factor'' $u(t)$ is appropriately accounted for.
By (\ref{ult}) Eq.~(\ref{fprod}), indeed, is equivalent to (\ref{convtheo}).

\subsection{Testing LT-consistency}
\label{SStesting}For a function $f(t)$ to be LT-consistent it is with regard to (\ref{varphi})
crucial that
\begin{equation}
u(0)f(0)=f(0).
\label{consunull}
\end{equation}
In this equation neither $u(0)$ nor $f(0)$ necessarily denote unique
definite values. The symbol $u(0)$
is just a synonym for the pseudo-function $u_0(t)$. As $f(t)$ may be
any linear combination of d-functions, causal functions, and derivatives,
$f(0)$ may actually denote any linear combination of definite unique values,
connect functions, and derivatives of the latter, i.e., impulse functions.
Keeping this in mind,
and observing the unit-step redundancy theorem (\ref{ustepequiv}),
the criterion for LT-consistency of any type of real function $f(t)$ can be
universally expressed by the

\parskip12pt\noindent
{\bf LT-consistency theorem:} The function $f(t)$ is LT-consistent if, and
only if,
\begin{equation}
u(t)f(t)=f(t)\hbox{\ \ for\ }-\infty<t<+\infty. 
\label{consu}
\end{equation}

The criterion (\ref{consu}) provides for the formal proof of the conclusion
that was already drawn from Fig.~\FigFunc, namely, that d-functions
$f_{\rm d}(t)$ are non-LT-consistent: $u(t)f_{\rm d}(t)\ne f_{\rm d}(t)$. 
By contrast, causal functions $f_{\rm c}(t)$ are LT-consistent, because
by the unit-step redundancy theorem there holds
$u(t)f_{\rm c}(t)=f_{\rm c}(t)$.

A useful application of the LT-consistency theorem is verification of 
LT-consistency of shifted functions.
From the formula for inverse L-transformation (\ref{ilt}) one obtains in the
familiar way the relationship
\begin{equation}
\varphi(t-\tau)=L^{-1}\{L\{f(t)\}\cdot{\rm e}^{-s\tau}\},
\label{varphitau}
\end{equation}
which indicates that multiplication of $L\{f(t)\}$ by $\exp(-s\tau)$ shifts
the causal inverse transform by any positive or negative amount $\tau$.
Utilizing (\ref{varphi})
one obtains from (\ref{varphitau}) by L-transformation the

\parskip8pt\noindent
{\bf shifting theorem}
\begin{equation}
L\{f(t)\}\cdot{\rm e}^{-s\tau}=L\{u(t-\tau)f(t-\tau)\};
\hbox{\ \ }-\infty<\tau<+\infty.
\label{theotrans}
\end{equation}
This expression holds for any type of real function $f(t)$ and for
any $\tau$. However, it is only for $\tau\ge0$ that (\ref{theotrans})
is consistent with the \textit{alias theorem}, indicating that
(\ref{theotrans}) is only for $\tau\ge0$ LT-consistent. Indeed,
for a shifted \textit{causal} function $f_{\rm c}(t-\tau)$
the criterion (\ref{consu}) reads
\begin{equation}
u(t)\cdot f_{\rm c}(t-\tau)=f_{\rm c}(t-\tau),
\label{consushift}
\end{equation}
and this condition is met only for $\tau\ge0$.

The shifting theorem makes particularly apparent that Laplace transformation
virtually is confined to causal functions. TLT's attempt to express the
shifting theorem in terms of ordinary derivable functions is awkward and
incoherent \cite{Doetsch1970,Doetsch1974}.

Yet, Laplace transformation of non-causal functions,
in particular, d-functions, does not entirely have to be ruled out. 
L-transformation of d-functions can be sensible and useful; cf.\
Sects.~\ref{SSdtsecond}, \ref{SSitsecond}, \ref{SSspontan}. One just has
to keep in mind that the L-transform of a d-function, $L\{f_{\rm d}(t)\}$,
actually is the L-transform of $u(t)f_{\rm d}(t)$, i.e.,
$L\{u(t)f_{\rm d}(t)\}$.

\subsection{Summary}
\label{SSgetsum}Laplace transformation of a real function $f(t)$ is said to be
consistent if the pertinent inverse L-transform is identical to $f(t)$.
For $t>0$ LT is in this sense consistent
for any type of real function that can be L-transformed at all.
However, such confined  consistency is not sufficient
for LT to provide a coherent system of operator calculus. Even in
the realm of TLT the behavior of functions at $t=0$ turns out in effect
to be involved. Explicit inclusion of the
point $t=0$ into LT's definition interval requires

\parskip8pt\noindent
a) appreciation
of the implicit bilaterality of the L-transformation formula (\ref{ult}); and

\parskip8pt\noindent
b) mathematical description of the behavior at $t=0$ of both $f(t)$ and
$\varphi(t)$, i.e., by the connect function $u_0(t)$ and/or its derivatives.

From a) there emerges a bilateral equivalent of the L-transformation
formula (\ref{ult}) that makes LT compatible with the theory of
Fourier-transformation.
This formula in turn is dependent on consistent definition of the
unit step function, i.e., according to b). L-transformation turns out to be
ambivalent, i.e., there holds $L\{f(t)\}=L\{u(t)f(t)\}$. Inverse
L-transformation is unequivocal; the inverse transform has
the form $u(t)f(t)$. Therefore, only causal functions are LT-consistent.

\section{The derivation/integration theorem}
\label{Stheo}There are two fundamentally different approaches to obtaining LT-theorems such
as those for derivation/integration:

\parskip8pt
\noindent
a) Determination by (\ref{ult}) of the L-domain operation that
corresponds to de\-ri\-va\-tion/in\-te\-gra\-tion of the {\em original
function} $f(t)$.

\parskip8pt
\noindent
b) Determination by (\ref{ilt}) of the L-domain operation that
corresponds to de\-ri\-va\-tion/in\-te\-gra\-tion of the
{\em inverse L-transform} $\varphi(t)$.

L-transformation by (\ref{ult}) or (\ref{bult}) is ambivalent, whereas
inverse L-transformation by (\ref{ilt}) is unequivocal. Moreover,
any kind of application of LT, in particular, to the solution of linear DEs,
is in the first place dependent on the inverse L-transform $\varphi(t)$
-- as opposed to the original function $f(t)$.
Therefore, only the approach b) is adequate.

TLT's derivation theorem (\ref{dttlt}) is based on the inadequate
approach a). In TLT,
the derivative $f'(t)$ of $f(t)$ is presupposed to exist in the ordinary
mathematical sense and the
L-transform of $f'(t)$ is expressed by (\ref{ult}). The theorem for the
first derivative then emerges from integration by parts:
\begin{eqnarray}
L\{f'(t)\}&=&\int_0^\infty f'(t){\rm e}^{-st}{\ \rm d}t\nonumber\\
&=&\Bigl[f(t){\rm e}^{-st}\Bigr]_0^\infty
+s\int_0^\infty f(t){\rm e}^{-st}{\ \rm d}t\nonumber\\
&=&sL\{f(t)\}-f(0).
\label{fdone}
\end{eqnarray}
Although this kind of mathematical reasoning is formally correct, the
result is not LT-consistent. The theorem
holds only for d-functions $f(t)=f_{\rm d}(t)$, i.e., functions that at
$t=0$ are derivable in the ordinary sense; d-functions are
non-LT-consistent because $u(t)f_{\rm d}(t)\ne f_{\rm d}(t)$; cf.\
(\ref{consu}).

Below, the LT-consistent theorems for derivation and integration are
obtained by the approach b). These theorems are termed the {\em primary}
derivation and integration theorem, respectively. They hold for
derivatives/integrals of the form $[u(t)f(t)]^{(\pm n)}$.
There also exist {\em secondary} theorems; these hold for $f^{(\pm n)}(t)$.
TLT's derivation theorem turns out to be of the secondary type. (In TLT
the distinction between $f^{(\pm n)}(t)$ and $[u(t)f(t)]^{(\pm n)}$
is ignored.) The primary theorems for derivation/integration have the same
form such that they can be unified. The unified theorem can be generalized
for non-integer order of derivation/integration. The
secondary theorems turn out to be virtually irrelevant.

\subsection{The primary derivation theorem}
\label{SSdtfirst}As the inverse transform $\varphi(t)$ is a causal and non-ordinary
function, its
derivatives are non-ordinary and causal, as well.
By contrast, the formula (\ref{ilt}) can within the integral be
derived for $t$ in the ordinary sense, and the order of derivation
is unlimited, as only the function $\exp(st)$ needs to be derived.
For the first derivative of $\varphi(t)$ one obtains
\begin{equation}
\varphi'(t)=L^{-1}\{sF(s)\}=L^{-1}\{sL\{f(t)\}\}=L^{-1}\{sL\{u(t)f(t)\}\}.
\label{pdtone}
\end{equation}
Observing the concatenation theorem (\ref{concat}), one obtains
for the second derivative 
\begin{equation}
\varphi''(t)=L^{-1}\{sL\{\varphi'(t)\}\}
=L^{-1}\{s^2L\{\varphi(t)\}\}
=L^{-1}\{s^2L\{f(t)\}\},
\label{pdttwo}
\end{equation}
Continuing with this kind of reasoning
to obtain higher-order derivatives, and observing
(\ref{varphi}), one obtains the

\noindent
{\bf primary derivation theorem}
\begin{equation}
L\{[u(t)f(t)]^{(n)}\}=s^nL\{f(t)\};\hbox{\ \ }n=0,1,\ldots
\label{dtfirst}
\end{equation}
Equation (\ref{dtfirst}) holds for any type of real function $f(t)$ and 
for any $n$ for which the inverse transform of $s^nL\{f(t)\}$ exists.

For causal functions $f_{\rm c}(t)$ the primary derivation theorem assumes
the form
\begin{equation}
L\{f_{\rm c}^{(n)}(t)\}=s^nL\{f_{\rm c}(t)\};\hbox{\ \ }n=0,1,\ldots
\label{dtfirstfc}
\end{equation}
The form (\ref{dtfirstfc}) applies, in particular, to the impulse functions
$f_{\rm c}(t)=\delta^{(n)}(t)$, and therefore
Eqs.~(\ref{implt}, \ref{deltanlt}) are consistent with (\ref{dtfirstfc}).
In contrast to
a widespread misconception, the reason why the form (\ref{dtfirstfc})
holds for impulse functions is not because these are distributions
but because they are causal.

The $t$-domain implications of the operation $s^nL\{f(t)\}$
are quite different for causal functions {\em versus} 
d-functions. As for causal functions $f_{\rm c}(t)$ there holds
$\varphi(t)=f_{\rm c}(t)$, there follows from the above deduction of
(\ref{dtfirst})
\begin{equation}
L^{-1}\{s^nL\{f_{\rm c}(t)\}\}=f_{\rm c}^{(n)}(t).
\label{fcninverse}
\end{equation}

For d-functions the corresponding relationship 
is considerably more complicated. To demonstrate this, functions of the type
\begin{equation}
f_{\rm ud}(t)=u(t)f_{\rm d}(t)
\label{fud}
\end{equation}
are taken into consideration, where $f_{\rm d}(t)$ denotes a d-function.
The $n$-th derivative of the so-called ud-function $f_{\rm ud}(t)$
is depicted by the

\parskip8pt
\noindent
{\bf dud-theorem}
\begin{equation}
f_{\rm ud}^{(n)}(t)=[u(t)f_{\rm d}(t)]^{(n)}
=u(t)f_{\rm d}^{(n)}(t)
+\sum_{\nu=0}^{n-1}f_{\rm d}^{(n-1-\nu)}(0)\cdot\delta^{(\nu)}(t);
\hbox{\ \ }n=1,2,\ldots
\label{udtheo}
\end{equation}
The dud-theorem ({\bf d}erivation of {\bf ud}-function) expresses the $n$-th
non-ordinary derivative of $u(t)f_{\rm d}(t)$
by the ordinary derivatives of $f_{\rm d}(t)$, i.e., at the expense of
getting impulse functions involved as depicted by (\ref{udtheo}). The
dud-theorem is explained in Sect.~\ref{SSudtheo}. 

By (\ref{udtheo}), utilizing (\ref{dtfirst}), the effect of the operation
$s^nL\{f_{\rm d}(t)\}$ gets depicted by
\begin{equation}
L^{-1}\{s^nL\{f_{\rm d}(t)\}\}=u(t)f_{\rm d}^{(n)}(t)
+\sum_{\nu=0}^{n-1}f_{\rm d}^{(n-1-\nu)}(0)\cdot\delta^{(\nu)}(t);
\hbox{\ \ }n=1,2,\ldots
\label{dtfirstfd}
\end{equation}

As an example, consider the operation $sL\{\cos\omega t\}$. From
(\ref{dtfirstfd}) one obtains
\begin{equation}
L^{-1}\{sL\{\cos\omega t\}\}=-u(t)\cdot\omega\sin\omega t+\delta(t).
\label{xcosa}
\end{equation}
The same result is obtained from
the LT-correspondences (\ref{corrsin}, \ref{corrcos}), i.e.,
\begin{equation}
L^{-1}\{sL\{\cos\omega t\}\}
=L^{-1}\Bigl\lbrace s\cdot{s\over s^2+\omega^2}\Bigr\rbrace
=L^{-1}\Bigl\lbrace 1-{\omega^2\over s^2+\omega^2}\Bigr\rbrace
=\delta(t)-u(t)\cdot\omega\sin\omega t.
\label{xcos}
\end{equation}

The most important domain of application of the primary derivation theorem is 
determination of the evoked solution of the inhomogeneous linear DE;
cf.\ Sect.~\ref{SSevoked}.

\subsection{The secondary derivation theorem}
\label{SSdtsecond}Although d-functions and their derivatives
are non-LT-consistent, the
question is legitimate how the L-transform of $f_{\rm d}^{(n)}(t)$
can be expressed by the L-transform of $f_{\rm d}(t)$. The
answer is already implied in (\ref{dtfirstfd}). From that expression one
obtains by L-transformation the

\parskip8pt
\noindent
{\bf secondary derivation theorem}
\begin{equation}
L\{f_{\rm d}^{(n)}(t)\}
=s^nL\{f_{\rm d}(t)\}
-\sum_{\nu=0}^{n-1}f_{\rm d}^{(n-1-\nu)}(0)\cdot s^{\nu};
\hbox{\ \ }n=1,2,\ldots
\label{dtsecond}
\end{equation}

The secondary derivation theorem (\ref{dtsecond}) turns out to be
identical to TLT's derivation theorem.
This is a consequence of TLT's original endeavour to
provide an operator calculus essentially for d-functions;
cf.\ Eq.~(\ref{fdone}).

The secondary derivation theorm is not a self-contained derivation
theorem as it just emerges from application of the primary derivation
theorem to d-functions.

\subsection{The primary integration theorem}
\label{SSitfirst}According to (\ref{ilt}) the inverse transform's integral function
can be depicted by
\begin{equation}
\varphi^{(-1)}(t)=L^{-1}\{s^{-1}F(s)\}=L^{-1}\{s^{-1}L\{f(t)\}\}
=L^{-1}\{s^{-1}L\{u(t)f(t)\}\}.
\label{pitone}
\end{equation}
As $\varphi(t)$ is causal and LT-consistent,
the integral function $\varphi^{(-1)}(t)$ is causal and LT-consistent, as
well. Thus for the second-order integral there holds
\begin{equation}
\varphi^{(-2)}(t)=L^{-1}\{s^{-1}L\{\varphi^{(-1)}(t)\}\}
=L^{-1}\{s^{-2}\{L\{\varphi(t)\}\}=L^{-1}\{s^{-2}L\{f(t)\}\};
\label{pittwo}
\end{equation}
and so on for higher-order integrals. Utilizing (\ref{varphi})
one obtains the

\parskip8pt
\noindent
{\bf primary integration theorem}
\begin{equation}
L\{[u(t)f(t)]^{(-n)}\}=s^{-n}L\{f(t)\};\hbox{\ \ }n=0,1,\ldots
\label{itfirst}
\end{equation}
which for causal functions $f(t)=f_{\rm c}(t)$ reads
\begin{equation}
L\{f_{\rm c}^{(-n)}(t)\}=s^{-n}L\{f_{\rm c}(t)\}.
\label{itfirstc}
\end{equation}

The operation $s^{-n}L\{f(t)\}$ is equivalent to $n$-fold integration of
the causal function $[u(t)f(t)]$, which implies $(n-1)$-fold iteration
of the $t$-domain operation
\begin{equation}
[u(t)f(t)]^{(-1)}=\int_{-\infty}^tu(\tau)f(\tau){\ \rm d}\tau
=u(t)\int_0^tf(\tau){\ \rm d}\tau.
\label{intoff}
\end{equation}
As $u(t)f(t)$ and the integral functions are causal, 
the $n$-fold integral can be expressed by the formula
\begin{equation}
[u(t)f(t)]^{(-n)}={u(t)\over(n-1)!}\int_0^t(t-\tau)^{n-1}f(\tau){\ \rm d}\tau;
\hbox{\ \ }n=1,2,\ldots
\label{nintform}
\end{equation}
The consistency of this formula with the primary integration theorem
(\ref{itfirst})
can be verified by L-transformation and application of the convolution
theorem (\ref{convtheo}), utilizing the LT-correspondence (\ref{corrtn}).

The $t$-domain implications of the operation $s^{-n}L\{f(t)\}$ are just
as different for causal functions \textit{versus} d-functions as was
found for derivation. For causal functions $f_{\rm c}(t)$ one obtains
\begin{equation}
L^{-1}\{s^{-n}L\{f_{\rm c}(t)\}\}=f_{\rm c}^{(-n)}(t).
\label{fcmn}
\end{equation}

When $f(t)=f_{\rm d}(t)$ is a d-function such that $f_{\rm d}(t)$ is
the $n$-th ordinary derivative of $f_{\rm d}^{(-n)}(t)$, there holds
\begin{eqnarray}
[u(t)f_{\rm d}(t)]^{(-1)}
&=&\int_{-\infty}^tu(t)f_{\rm d}(\tau){\ \rm d}\tau
=u(t)\int_0^tf_{\rm d}(\tau){\ \rm d}\tau\nonumber\\
&=&u(t)[f_{\rm d}^{(-1)}(t)-f_{\rm d}^{(-1)}(0)].
\label{udintone}
\end{eqnarray}
By iteration of (\ref{udintone}) one obtains the so-called

\parskip8pt
\noindent
{\bf iud-theorem} ({\bf i}ntegration of {\bf ud}-function)
\begin{equation}
f_{\rm ud}^{(-n)}(t)=[u(t)f_{\rm d}(t)]^{(-n)}=u(t)f_{\rm d}^{(-n)}(t)
-u(t)\sum_{\nu=0}^{n-1}f_{\rm d}^{(-n+\nu)}(0)\cdot{t^\nu\over\nu!};
\hbox{\ \ }n=1,2,\ldots
\label{iudtheo}
\end{equation}

For the effect of the operation $s^{-n}L\{f_{\rm d}(t)\}$
one obtains from (\ref{itfirst}) and (\ref{iudtheo})
\begin{equation}
L^{-1}\{s^{-n}L\{f_{\rm d}(t)\}\}=u(t)f_{\rm d}^{(-n)}(t)
-u(t)\sum_{\nu=0}^{n-1}f_{\rm d}^{(-n+\nu)}(0)\cdot{t^\nu\over\nu!};
\hbox{\ \ }n=1,2,\ldots
\label{iudtheolt}
\end{equation}

As an example, consider the operation $s^{-1}L\{\exp(-at)\}$.
From (\ref{iudtheolt}) one obtains
\begin{equation}
L^{-1}\{s^{-1}L\{{\rm e}^{-at}\}\}
=-{u(t)\over a}\cdot{\rm e}^{-at}+{u(t)\over a}.
\label{xexpata}
\end{equation}
Using the 
LT-cor\-res\-pon\-den\-ces (\ref{correxp}, \ref{correxpi}) one obtains
the same result, i.e.,
\begin{equation}
L^{-1}\{s^{-1}L\{{\rm e}^{-at}\}\}
=L^{-1}\Bigl\lbrace{1\over s}\cdot{1\over s+a}\Bigr\rbrace
={u(t)\over a}(1-{\rm e}^{-at}).
\label{xexpat}
\end{equation}

\subsection{The secondary integration theorem}
\label{SSitsecond}By analogy to the secondary derivation theorem, Eq.~(\ref{iudtheolt})
enables for  d-functions expression of the
L-transform of the $n$-th order integral in terms of the L-transform
of the d-function itself. By L-transformation of (\ref{iudtheolt}),
utilizing (\ref{corrtn}), there emerges the

\parskip8pt
\noindent
{\bf secondary integration theorem}
\begin{equation}
L\{f_{\rm d}^{(-n)}(t)\}
=s^{-n}L\{f_{\rm d}(t)\}
+\sum_{\nu=0}^{n-1}f_{\rm d}^{(-n+\nu)}(0)\cdot s^{-\nu-1};
\hbox{\ \ }n=1,2,\ldots
\label{itsecond}
\end{equation}
The secondary integration theorem is complementary to the secondary derivation
theorem. For $n=1$ the two theorems are equivalent; indeed,
from (\ref{itsecond}) one obtains
\begin{equation}
L\{f_{\rm d}^{(-1)}(t)\}=s^{-1}L\{f_{\rm d}(t)\}
+f_{\rm d}^{(-1)}(0)\cdot s^{-1}.
\label{iudti}
\end{equation}
Taking into account that by definition $f_{\rm d}(t)$ is the first ordinary
derivative of $f_{\rm d}^{(-1)}(t)$, 
Eq.~(\ref{iudti}) is identical to the secondary derivation theorem and thus
to TLT's derivation theorem.

As TLT's derivation theorem (\ref{dttlt}, \ref{fdone}) is identical
to the secondary
derivation theorem (\ref{dtsecond}), one would expect TLT's integration theorem
to comply with the secondary integration theorem (\ref{itsecond}). However,
TLT's integration theorem (\ref{ittlt}) actually is identical to the
{\em primary} integration theorem (\ref{itfirst}). The latter
theorem holds only for causal functions; cf.\ Eq.~(\ref{itfirstc});
its application
to d-functions -- such as is customary in TLT -- will in general yield
erroneous results.

Thus, it turns out that the alleged mathematical consistency of TLT's
theorems for derivation/integration is delusive. 
TLT's derivation theorem holds only for d-functions whereas
TLT's integration theorem holds only for causal functions. In TLT both
theorems are ordinarily used in an untenable way: TLT's derivation theorem
is regarded as \textit{the} derivation theorem of LT although it
holds only for the non-LT-consistent d-functions.
TLT's integration theorem is ordinarily utilized
for d-functions although it does not apply to this type of function.
Eventually, these observations explain why there is a formal conflict
between TLT's theorems for derivation and integration
(\ref{dttlt}, \ref{ittlt}).

The secondary integration theorem is not a self-contained integration theorem
as it just emerges from application of the primary integration theorem to
d-functions.

\subsection{The generalized derivation/integration theorem}
\label{SSgditheo}The LT-consistent theorems for derivation and integration, i.e.,
(\ref{dtfirst}) and (\ref{itfirst}), can obviously be unified into one
formula, i.e.,
\begin{equation}
L\{[u(t)f(t)]^{(n)}\}=s^{n}L\{f(t)\}\hbox{\ \ }n=0,\pm1,\pm2,\ldots
\label{ditheo}
\end{equation}
Thus, ultimately $t$-domain integration is in the L-domain merely 
``reciprocal'' to $t$-domain derivation. The $t$-domain implications
of the operation $s^{n}L\{f(t)\}$ that are outlined in
Sects.~\ref{SSdtfirst}, \ref{SSitfirst} must be observed.

As a consequence of the formal equivalence of derivation and integration in
the L-domain description, the theorem can be generalized for
non-integer order of
derivation/integration. For any real $r\ge0$ there holds
(Riemann, Liouville, Cauchy)
\begin{equation}
D^{-r}\{u(t)f(t)\}={u(t)\over\Gamma(r)}\int_0^t
(t-\tau)^{r-1}f(\tau){\ \rm d}\tau;\hbox{\ }r\ge0.
\label{rint}
\end{equation}
The operator $D$ denotes generalized derivation/integration,
while the negative exponent $-r$ indicates $r$-th-order {\em integration}.
As the integral in (\ref{rint}) is equivalent to the
convolution $t^{r-1}\ast f(t)$, Eq.~(\ref{rint}) can
by the convolution theorem (\ref{convtheo}) be L-transformed,
and when the LT-correspondence (\ref{corrtr}) is utilized one obtains
\begin{eqnarray}
L\{D^{-r}\{u(t)f(t)\}\}&=&
{1\over\Gamma(r)}L\{t^{r-1}\}\cdot L\{f(t)\}\\
&=&s^{-r}L\{f(t)\};\hbox{\ }r\ge0.
\label{rinttheo}
\end{eqnarray}

Equation (\ref{rinttheo}) depicts the generalized {\em integration} theorem.
This theorem suffices to obtain the same kind of
generalization for derivation, namely, by
concatenation of (\ref{rinttheo}) and (\ref{dtfirst}).
The sequence of $r$-th order integration and
$n$-th order derivation yields $(n-r)$th order derivation or integration,
depending on whether $n>r$ or $n<r$ \cite{Sokolov}.
Denoting $n-r=\alpha$ one obtains the

\parskip8pt\noindent
{\bf gdi-theorem} ({\bf g}eneralized {\bf d}erivation/{\bf i}ntegration theorem)
\begin{equation}
L\{D^{\alpha}\{u(t)f(t)\}\}=s^\alpha L\{f(t)\};\hbox{\ \ }\alpha\in\Re. 
\label{gditheo}
\end{equation}
In the form (\ref{gditheo}) the theorem holds for any type of real function
$f(t)$.
For integer values of $\alpha$ Eq.~(\ref{gditheo}) is equivalent to
(\ref{ditheo}).

The concept of \textit{fractional calculus}, i.e., utilization of
derivatives of non-integer order, is essentially based on the
formula (\ref{rint}), and this formula holds only for causal functions.
In the realm of TLT derivatives of non-integer order
can not be consistently expressed at all, as
TLT's derivation theorem neither holds for causal functions nor for
non-integer order of derivation. When TLT's derivation theorem is used anyway, 
initial values become involved which merely are a nuisance.
It is by the present approach
-- which provides for the generalized theorem (\ref{ditheo}) -- that
Laplace transformation becomes a consistent and efficient tool for doing
fractional calculus.

Examples for the application of (\ref{gditheo}), i.e., for $\alpha=1/2$,
are listed in Sect.~\ref{SSxdhalf}.

\subsection{Summary}
\label{SSdtsum}LT-consistent theorems for t-domain derivation and integration are
obtained by deduction from the behavior of the inverse L-transform.
These theorems -- the \textit{primary} theorems -- are formally congruent
and thus can be unified. The unified primary theorem is
generalized for non-integer order of derivation/integration.
The generalized theorem (the gdi-theorem) accounts in the
L-domain for both integer-order and non-integer-order of $t$-domain
derivation/integration. As the inverse transform is causal
the primary theorems and the gdi-theorem hold for causal functions.

The respective secondary theorems emerge from application of the primary ones
to the L-transforms of ordinary derivable functions (d-functions).
TLT's derivation theorem is identical to the secondary
derivation theorem. This theorem can not be regarded as \textit{the}
derivation theorem of Laplace transformation because it applies only to the
class of d-functions, which is non-LT-consistent. TLT's integration theorem is
identical to
the primary integration theorem and therefore holds only for causal
functions; its utilization for d-functions -- such as in TLT -- is untenable.

\section{Solution of linear differential equations}
\label{Ssolu}Below, the new approach to the solution of linear DEs by LT
is demonstrated for a fairly general form of the ordinary linear DE, namely,
\begin{equation}
\sum_{n=0}^Na_{\rm n}y^{(n)}(t)=\sum_{m=0}^Mb_{\rm m}x^{(m)}(t);
\hbox{\ \ \ }N=1,2\ldots;\hbox{\ }M=0,1\ldots
\label{diffeq}
\end{equation}
The coefficients $a_{\rm n}$, $b_{\rm m}$ are presupposed to be real
constants. Two functions are involved, namely, the
\textit{excitation} function $x(t)$ and the \textit{response function}
$y(t)$. The form
(\ref{diffeq}) is more general than that accounted for by the TLT method,
as in (\ref{diffeq}) derivatives of the excitation function are allowed
to be included, i.e., corresponding to $M>0$.

A clear distinction is made between the
\textit{evoked} solution and the \textit{spontaneous} solution of
the inhomogeneous DE \cite{Terhardt1987}. The evoked solution (the
pertinent system's evoked response to $x(t)$) is that part of the total
response
that is elicited by the excitation function $x(t)$ alone. The spontaneous
solution (the system's ``spontaneous'' response), if it exists, is
ascribed to a pre-excited initial state of the system. Mathematically,
the evoked solution is the DE's total (general) solution minus the
general solution of the pertinent homogeneous DE. The spontaneous solution
is identical and synonymous to the general solution of the homogeneous DE.

\subsection{The evoked response:\\
Particular solution of the inhomogeneous DE}
\label{SSevoked}For the excitation function $x(t)$ to be LT-consistent it must be
defined as a causal function. As a consequence, its derivatives have
the form $[u(t)x(t)]^{(m)}$. Another consequence is that
the evoked response $y_{\rm e}(t)$ is causal as well; thus, its
derivatives implicitly assume the form $[u(t)y_{\rm e}(t)]^{(n)}$.
Therefore, conversion of (\ref{diffeq}) into the L-domain is governed by
the primary derivation theorem (\ref{dtfirst}). 
The L-domain representation of (\ref{diffeq}) reads
\begin{equation}
\sum_{n=0}^Na_{\rm n}s^nL\{y_{\rm e}(t)\}
=\sum_{m=0}^Mb_{\rm m}s^mL\{x(t)\};
\hbox{\ \ }N=1,2\ldots;\hbox{\ }M=0,1\ldots,
\label{diffeqult}
\end{equation}
and one eventually obtains the evoked solution
\begin{equation}
y_{\rm e}(t)=L^{-1}\Biggl\lbrace
{\sum_{m=0}^Mb_{\rm m}s^m\over\sum_{n=0}^Na_{\rm n}s^n}\cdot L\{x(t)\}
\Biggr\rbrace.
\label{diffeqltsol}
\end{equation}

In contrast to the TLT method, the evoked response
is by (\ref{diffeqltsol}) obtained without involvement of initial values, i.e.,
additional constants. To
obtain the evoked response one does not have to pretend that the pertinent
system is in a non-preexcited initial state. With respect to the evoked
response the system's initial state is irrelevant.

For the operation $L^{-1}$, i.e., inverse L-transformation, there exist two
alternatives. The first of them is taking (\ref{diffeqltsol}) at face value,
which implies that the L-transform of $x(t)$ has to be included. The second
alternative exploits the convolution theorem (\ref{convtheo}), whereby
inverse L-transformation can be confined to the first factor in
(\ref{diffeqltsol}):
\begin{equation}
h(t)=L^{-1}\Biggl\lbrace
{\sum_{m=0}^Mb_{\rm m}s^m\over\sum_{n=0}^Na_{\rm n}s^n}\Biggr\rbrace.
\label{impresp}
\end{equation}
The function $h(t)$ depicts the evoked response to the excitation function
$x(t)=\delta(t)$, i.e., the impulse response.
As (\ref{diffeqltsol}) has the form
\begin{equation}
y_{\rm e}(t)=L^{-1}\{L\{h(t)\}\cdot L\{x(t)\}\},
\label{desolform}
\end{equation}
by (\ref{convtheo}) there emerges the familiar convolution formula
\begin{equation}
y_{\rm e}(t)=u(t)\int_0^th(t-\tau)x(\tau){\ \rm d}\tau.
\label{convint}
\end{equation}

For example, for $N=1$, $M=0$ one obtains the impulse response
\begin{equation}
h(t)=L^{-1}\Biggl\lbrace{b_0\over a_0+a_1s}\Biggr\rbrace
=u(t)\beta{\rm e}^{-\alpha t};\hbox{\ \ }\alpha=a_0/a_1;\hbox{\ }\beta=b_0/a_1.
\label{imprespfo}
\end{equation}
Utilizing (\ref{imprespfo}) and (\ref{convint}), one obtains the

\parskip=8pt\noindent
{\bf evoked solution of the first-order linear DE} 
\begin{equation}
y_{\rm e}(t)=u(t)\beta\int_0^tx(\tau){\rm e}^{-\alpha(t-\tau)}{\ \rm d}\tau
=u(t)\beta{\rm e}^{-\alpha t}\int_0^tx(\tau){\rm e}^{\alpha\tau}{\ \rm d}\tau.
\label{evoksolfo}
\end{equation}
Though this solution is entirely expressed in the $t$-domain, it
nevertheless is based on LT, which implies that (\ref{evoksolfo}) is correct
only for causal excitation functions $x(t)$.

It may be noted that it is the independent LT-based expression
of the evoked solution (\ref{diffeqltsol}) that provides a solid basis to
\textit{operator calculus} as it is customarily employed
in the theories of linear control systems and electrical circuits.

\subsection{The spontaneous response:\\
General solution of the homogeneous DE}
\label{SSspontan}When the evoked response is by LT determined as just described,
for the solution of the homogeneous DE still any method available can be
chosen. From the present approach there emerges a new LT-based method which
exploits the facts that

\parskip8pt
\noindent
a) the spontaneous response $y_{\rm s}(t)$
and its derivatives are d-functions; their L-transforms have the form
$L\{y_{\rm s}^{(n)}(t)\}=L\{u(t)y_{\rm s}^{(n)}(t)\}$;

\parskip8pt
\noindent
b) the homogeneous  DE
can by the dud-theorem be converted into an equivalent inhomogeneous DE.

To get the homogeneous DE made up for conversion into an equivalent
inhomogeneous DE, it is multiplied by $u(t)$; this yields
\begin{equation}
a_0u(t)y_{\rm s}(t)+a_1u(t)y_{\rm s}'(t)+\ldots+a_Nu(t)y_{\rm s}^{(N)}(t)=0.
\label{dehomou}
\end{equation}
Because of LT's ambivalence,
the L-transform of (\ref{dehomou}) is identical to the L-transform of the
original homogeneous DE. Application of the dud-theorem (\ref{udtheo})
converts (\ref{dehomou}) into the form
\begin{eqnarray}
a_0u(t)y_{\rm s}(t)&+&
a_1\Bigl\langle[u(t)y_{\rm s}(t)]'-y_{\rm s}(0)\cdot\delta(t)\Bigr\rangle
+a_2\Bigl\langle[u(t)y_{\rm s}(t)]''-y_{\rm s}'(0)\cdot\delta(t)
-y_{\rm s}(0)\cdot\delta'(t)\Bigr\rangle\nonumber\\
&+&\ldots
+a_N\Bigl\langle[u(t)y_{\rm s}(t)]^{(N)}
-\sum_{\nu=0}^{N-1}y_{\rm s}^{(N-1-\nu)}(0)\cdot\delta^{(\nu)}(t)\Bigr\rangle=0.
\label{dehomoud}
\end{eqnarray}
Equation (\ref{dehomoud}) is equivalent to (\ref{dehomou}). However,
(\ref{dehomoud}) actually is
an {\em inhomogeneous} DE. The impulse functions included
in (\ref{dehomoud}) play the role of virtual excitation functions
\cite{Terhardt1987}.
This becomes particularly apparent when (\ref{dehomoud}) is expressed
in the form
\begin{equation}
\sum_{n=0}^Na_n[u(t)y_{\rm s}(t)]^{(n)}
=\sum_{\mu=0}^{N-1}c_\mu\delta^{(\mu)}(t),
\label{dehomoue}
\end{equation}
where the coefficients $c_\mu$ are depicted by
\begin{equation}
c_\mu=\sum_{\nu=0}^{N-1-\mu}a_{\mu+\nu+1}\cdot y_{\rm s}^{(\nu)}(0).
\label{cmu}
\end{equation}
From (\ref{dehomoue}) the solution of the homogeneous DE can be obtained
as an {\em evoked} solution, i.e., as described in Sect.~\ref{SSevoked}.
One eventually obtains
\begin{equation}
u(t)y_{\rm s}(t)=L^{-1}\Biggl\lbrace
{\sum_{\mu=0}^{N-1}c_\mu s^{\mu}\over\sum_{n=0}^Na_ns^n}\Biggr\rbrace.
\label{sponsollt}
\end{equation}

The notation $u(t)y_{\rm s}(t)$ is not redundant, because $y_{\rm s}(t)$ is a
d-function. Equation (\ref{sponsollt}) in fact depicts
the {\em evoked} response of a \textit{virtual system}, i.e., to the
excitation function $\delta(t)$. The transmission function of that
virtual system is defined by the quotient in (\ref{sponsollt}), whose
numerator essentially originates from the impulse functions contained in
(\ref{dehomoue}).

Finally, the spontaneous solution $y_{\rm s}(t)$ itself can
by extrapolation into $t\le0$
be obtained from (\ref{sponsollt}). As according to the causality
theorem (\ref{varphi}) the operation $L^{-1}$ on the right side of
(\ref{sponsollt}) yields either explicitly or
implicitly a function of the form $u(t)y_{\rm s}(t)$, that kind of
extrapolation is equivalent to cancellation
of $u(t)$ on both sides of (\ref{sponsollt}), i.e., after inverse
L-transformation.

For example, for $N=1$ one obtains from (\ref{sponsollt}) and (\ref{cmu})
\begin{equation}
u(t)y_{\rm s}(t)
=L^{-1}\Biggl\lbrace{a_1y_{\rm s}(0)\over a_0+a_1s}\Biggr\rbrace
=u(t)y_{\rm s}(0){\rm e}^{-\alpha t};\hbox{\ \ }\alpha=a_0/a_1.
\label{usponsolfo}
\end{equation}
By cancellation of $u(t)$ there emerges from (\ref{usponsolfo}) the

\parskip=8pt\noindent
{\bf spontaneous solution of the first-order linear DE}
\begin{equation}
y_{\rm s}(t)=y_{\rm s}(0){\rm e}^{-\alpha t}\hbox{\ \ for\ }-\infty<t.
\label{sponsolfo}
\end{equation}

The described deduction of (\ref{sponsollt}) ultimately constitutes
a new method for obtaining the general solution of
the linear homogeneous DE of finite order $N$, i.e., by Laplace
transformation. The method
involves a) determination by (\ref{cmu})
of the coefficients $c_\mu$, i.e., from the DE's coefficients $a_n$; b)
insertion of the coefficients $c_\mu$ into (\ref{sponsollt}); c)
inverse L-transformation; and d) cancellation of the factor $u(t)$.
The $N$ arbitrary constants which invariably get involved are provided by the
initial values at $t=0$ of the solution $y_{\rm s}(t)$ itself and of
the latter's $N-1$ derivatives. The method is straightforward; yet,
mathematical intricacies may occur
in step c), i.e., inverse L-transformation.

As a pragmatic alternative to the described deduction,
the formula (\ref{sponsollt}) may
be obtained in a more immediate way, namely, by using for LT-conversion of
the homogeneous DE the secondary derivation theorem (\ref{dtsecond}).
This option is enabled by the fact that
the secondary derivation theorem accounts for the combination of
the dud-theorem with the
primary derivation theorem; cf.\ Sect.~\ref{SSdtsecond}. As the secondary
derivation theorem is identical to TLT's derivation theorem, this option
explains why TLT's total solution of the inhomogeneous DE is quite similar
to the new total solution; see the following two sections.

\subsection{The total solution}
\label{SStotal}Assuming that for the solution of the homogeneous DE the method just
described is chosen, the total solution of (\ref{diffeq}) can be compactly
depicted by a formula, i.e.,
\begin{equation}
y(t)=L^{-1}\Biggl\lbrace
{\sum_{m=0}^Mb_{\rm m}s^m\over\sum_{n=0}^Na_{\rm n}s^n}\cdot
L\{x(t)\}\Biggr\rbrace
+L_{\rm d}^{-1}\Biggl\lbrace
{\sum_{\mu=1}^Na_\mu\sum_{\nu=0}^{\mu-1}y_{\rm s}^{(\mu-1-\nu)}(0)s^{\nu}\over
\sum_{n=0}^Na_{\rm n}s^n}\Biggr\rbrace.
\label{desolg}
\end{equation}
The first term on the right side of (\ref{desolg}) depicts the
evoked response and is identical to (\ref{diffeqltsol}). The second term
depicts the spontaneous response by one single expression; this term
is equivalent to the combination of Eqs.(\ref{sponsollt}) and
(\ref{cmu}). The operator $L^{-1}$ indicates inverse transformation by
(\ref{ilt}). The operator $L^{-1}_{\rm d}$ indicates inverse transformation
by (\ref{ilt}) followed by extrapolation into $t\le0$, i.e.,
cancellation of the factor $u(t)$.

When the inverse L-transform is looked up from a TLT-based
table of LT-cor\-res\-pon\-den\-ces
it must be observed that in those tables the $t$-domain functions
are depicted in
the non-causal form. In Sect.~\ref{SScorr} examples are listed
of the LT-consistent notation of $t$-domain functions.

It must be kept in mind that LT-based solutions of linear DEs hold
only for \textit{causal}
excitation functions $x(t)=u(t)x(t)$. The solution for the ``steady-state''
may be obtained by asymptotical approach, i.e., for $t\to\infty$.

An example for application of (\ref{desolg}) is depicted in
Sect.~\ref{SSxcirc}.

Finally, it should be noticed that
for $N=1$, $M=0$ the LT-based solution of (\ref{diffeq}) can be entirely
expressed in the $t$-domain. By superposition of Eqs.~(\ref{evoksolfo}) and
(\ref{sponsolfo}) one obtains the

\parskip=8pt\noindent
{\bf total solution of the first-order linear DE}
\begin{equation}
y(t)=y_{\rm e}(t)+y_{\rm s}(t)
=u(t)\beta{\rm e}^{-\alpha t}\int_0^t x(\tau){\rm e}^{\alpha\tau}{\ \rm d}\tau
+y_{\rm s}(0){\rm e}^{-\alpha t};
\hbox{\ \ }\alpha=a_0/a_1;\hbox{\ }\beta=b_0/a_1.
\label{totsolfo}
\end{equation}
Thus, for the special case that in
(\ref{diffeq}) there is $N=1$, $M=0$ one does not actually
have to do L-transformation at all. The solution is
reduced to evaluation of the integral contained in (\ref{totsolfo}) and it
holds for causal excitation functions $x(t)=u(t)x(t)$.

\subsection{TLT's initial-value conflict}
\label{SSthing}One of the most annoying deficiencies of TLT is potential interference
of the DE's evoked solution with the spontaneous solution. The danger of
this kind of interference to occur emerges from the fact that the
DE's solution
obtained by TLT includes an inappropriate type of initial values. This is
why the phenomenon is termed the initial-value conflict.

When the general solution of a linear DE of the form (\ref{diffeq})
is worked out by means of TLT -- which is possible for $M=0$ -- one
eventually obtains the formula
\begin{equation}
y(t)=L^{-1}\Biggl\lbrace
{b_0\over\sum_{n=0}^Na_{\rm n}s^n}\cdot L\{x(t)\}
+{\sum_{\mu=1}^Na_\mu\sum_{\nu=0}^{\mu-1}y^{(\mu-1-\nu)}(0)s^{\nu}\over
\sum_{n=0}^Na_{\rm n}s^n}\Biggr\rbrace.
\label{desolgtlt}
\end{equation}
This formula is quite similar to (\ref{desolg}), i.e., for $M=0$. (The
constant $b_0$ is in (\ref{desolgtlt}) included just for formal
compatibility with (\ref{desolg}); one may in both equations
assume $b_0=1$.) Letting
aside the difference in inverse L-transformation that is apparent by
comparison of (\ref{desolgtlt}) to (\ref{desolg}), there remains the
crucial difference that (\ref{desolg}) correctly contains
the initial values of the \textit{spontaneous} solution, $y_{\rm s}^{(\nu)}(0)$,
whereas (\ref{desolgtlt}) contains the initial values of the \textit{total}
solution, $y^{(\nu)}(0)$. The latter, inadequate kind of initial
values inevitably emerge from LT-conversion of the DE (\ref{diffeq}) by
TLT's derivation theorem.

The potential initial-value conflict arises from the fact that the
evoked and the spontaneous solutions are superimposed in $y(t)$, such that
\begin{equation}
y^{(\nu)}(0)=y_{\rm e}^{(\nu)}(0)+y_{\rm s}^{(\nu)}(0);
\hbox{\ \ }\nu=0,1,\ldots
\label{ynnull}
\end{equation}
The initial values of the evoked solution, $y_{\rm e}^{(\nu)}(0)$, depend on
$x(t)$ and on the DE's coefficients $a_{\rm n}$, as is depicted by 
(\ref{diffeqltsol}). If there happens to be $y_{\rm e}^{(\nu)}(0)\equiv0$, i.e.,
for all $\nu=0,1,\ldots,N-1$, then there holds
$y^{(\nu)}(0)\equiv y_{\rm s}^{(\nu)}(0)$,
and (\ref{desolgtlt}) will yield the same result as (\ref{desolg}),
i.e., except for the 
difference in inverse L-transformation, and for $M=0$.
However, for many types of system
and of excitation function $y_{\rm e}^{(\nu)}(0)$ may turn out to be
different from null, i.e., for
at least one value of $\nu$. If this occurs then the result obtained by TLT is
incorrect; the initial values $y^{(\nu)}(0)$ can no longer be
freely chosen, i.e., independent of the excitation function; cf.\
Sect.~\ref{SSfail}.

Obviously, (\ref{desolgtlt}) can easily be ``patched'', i.e.,
by arbitrarily substituting
the values $y_{\rm s}^{(\nu)}(0)$ for $y^{(\nu)}(0)$.
Most remarkably, there exists another, less trivial
kind of patch, namely, arbitrary substitution in (\ref{desolgtlt}) of
the values $y^{(\nu)}(-0)$ for $y^{(\nu)}(0)$.
The reason why the latter kind of patch works around the
initial-value conflict originates from the fact that the evoked response
and its derivatives are causal (cf.\ Sect.~\ref{SSevoked}),
such that there holds $y_{\rm e}^{(\nu)}(-0)\equiv0$, i.e., even when
$y_{\rm e}^{(\nu)}(0)\ne0$. As
the spontaneous response and its derivatives are d-functions, there holds
$y_{\rm s}^{(\nu)}(-0)\equiv y_{\rm s}^{(\nu)}(0)$. Thus, one obtains
from (\ref{ynnull})
\begin{equation}
y^{(\nu)}(-0)=y_{\rm s}^{(\nu)}(0)\hbox{\ \ for\ }\nu=0,1,\ldots,N-1.
\label{ynullminus}
\end{equation}
The identity (\ref{ynullminus}) holds for any type of linear DE and of
excitation function, From this identity there follows that
substitution in (\ref{desolgtlt}) of the initial point
$t=-0$ for $t=0$ renders (\ref{desolgtlt}) compatible with (\ref{desolg}),
such that TLT's solution becomes essentially correct,
i.e., except for the difference in inverse L-transformation.
The ``empirical'' finding that by this kind of patch TLT's 
initial-value conflict actually is worked around, becomes in this way
explained.

It should be noticed that the latter kind of patch is just as arbitrary as
the former. Even if one had no doubt at all about TLT's
consistency, the mere
fact that TLT's solution (\ref{desolgtlt}) in general needs to be patched
should have elicited that kind of doubt. One can not in earnest be
satisfied with either patch's working around TLT's initial-value conflict.

\subsection{Summary}
\label{SSsolsum}The total solution of the inhomogeneous linear DE is obtained by two
independent algorithms, i.e., one for the evoked response, another for the
spontaneous response.
As both the excitation function and the evoked response function are
causal, the derivatives of these functions
assume the form $[u(t)x(t)]^{(m)}$ and
$[u(t)y_{\rm e}(t)]^{(n)}$, respectively. When the inhomogeneous DE
is by the primary derivation theorem converted into the L-domain one
obtains a linear algebraic equation that does not contain any arbitrary
constants. This equation is 
solved in the familiar way to obtain the evoked solution
$y_{\rm e}(t)$ which is defined for $-\infty<t$.

For the determination of the {\em spontaneous response} $y_{\rm s}(t)$,
i.e., solution of the pertinent homogeneous DE, the fact that
the response function and its derivatives are d-functions is \textit{a priori}
taken into account. The L-transforms of these functions assume the form
$L\{y_{\rm s}^{(n)}(t)\}=L\{u(t)y_{\rm s}^{(n)}(t)\}$. As a consequence, the
homogeneous DE can be multiplied by $u(t)$ without affecting its L-transform.
The homogeneous DE such modified can by the dud-theorem
be converted into an equivalent inhomogeneous DE.
The L-transform of the latter DE depicts a virtual linear system
whose impulse response equals for $t>0$
the spontaneous solution of the original homogeneous DE.
By extrapolation into $t\le0$ one obtains from that impulse response
the spontaneous solution $y_{\rm s}(t)$ for $-\infty<t$.

For $M=0$, i.e., if the excitation function does not include any derivative,
the total solution obtained by TLT is formally similar to the correct solution 
depicted by the new method. However, the TLT-based solution suffers
from an initial-value conflict. This conflict can be worked around
by manipulation of the initial values.
The necessity of patching TLT's solution is another
symptom of TLT's deficiency.

\section{Concluding remarks}
\label{SSconc}From the theorems pointed out in Sect.~\ref{Sget} it should be apparent
that the key to the consistent LT-theory lies in the hidden
bilaterality of the transformation integral (\ref{ult}). By strictly
taking that bilaterality into account there emerges a structure
of LT-theory that is fundamentally different from that of TLT. Yet,
many of the new theory's theorems and expressions are familiar from TLT.
Pronounced formal differences from TLT's expressions occur, e.g.,
in the context of the theorems for derivation and integration.
The difference between the new formula (\ref{desolg}) for the general
solution of the linear inhomogeneous DE from the corresponding formula
(\ref{desolgtlt}) obtained by TLT, though at first sight marginal,
actually is crucial, namely, for resolution of TLT's inherent
initial-value conflict.

Mathematical consistency of the new approach to LT is achieved at the expense
of non-ordinary $t$-domain functions being crucially involved --
primarily, the inverse
L-transform $\varphi(t)$, which includes the connect function $u_0(t)$
and, possibly, its (non-ordinary) derivatives. This notion brings to mind the
aspect of mathematical rigour. If one maintains that
there is no rigorous mathematical account for non-ordinary functions and
non-ordinary derivatives then
the inverse L-transform can not in general be rigorously accounted for and
a rigorous theory of Laplace transformation can not exist at all.

Rigorous treatment of non-ordinary functions may be achieved by invoking the
theory of distributions. However, what makes LT-theory consistent
in the first place is adequate incorporation of non-ordinary functions and
non-ordinary derivatives. Mere attachment of distribution theory to TLT
can not cure TLT's inconsistency.
Mathematically rigorous treatment of non-ordinary functions and derivatives
may be achieved separately and subsequently -- which, historically, actually
has happened to the delta-impulse.

Making in this sense a distinction between mathematical consistency and
rigour, one can say that the new approach to Laplace transformation
outlined in the present article is consistent, even though it
may not fully qualify as mathematically rigorous. This
achievement is distinctly
preferable to inconsistency disguised by fake mathematical rigour --
which virtually is what is offered by Doetsch \cite{Doetsch1970,Doetsch1974}
and many others.

Once one has developed some scepticism about TLT one may be prepared to
realizing that
Doetsch's books are written in a somewhat defensive and dogmatic style,
and that Doetsch devoted considerable portions of text to
reasoning by plausibility and to explaining away apparent discrepancies.
From these observations it may be concluded that Doetsch virtually
was aware of TLT's inherent problems. Indeed, in the preface to
\cite{Doetsch1955} he suggested the need for a fundamental redesign of TLT.
However, until the end of his life (1977) he
stuck to the original layout of TLT, i.e., with some distribution theory as a
complement. In the two
volumes of his \textit{Handbook} \cite{Doetsch1950,Doetsch1955}, Doetsch
assembled a wealth of mathematics on LT -- which, however, may
distract from TLT's fundamental flaws.

To obtain a more thorough understanding of TLT's unfortunate history
it is probably
helpful to take notice of Doetsch's biography \cite{Connor2004,Remmert1999}.

\begin{appendix}
\section{APPENDIX}
\subsection{Failure of TLT: An example DE}
\label{SSfail}TLT promises to provide a method for obtaining the general solution of the
linear inhomogeneous DE, i.e., for any kind of excitation
function $x(t)$. Below, this promise is challenged by an example, i.e.,
for the first-order DE
\begin{equation}
ay(t)+y'(t)=\hat x\sin\omega t,
\label{xde}
\end{equation}
where $\hat x\sin\omega t$ denotes the excitation function and $y(t)$
denotes the response function.

From the theory of linear DE's one obtains the general solution
\begin{equation}
y(t)={\hat x\over a^2+\omega^2}(a\sin\omega t-\omega\cos\omega t)
+y_{\rm s}(0){\rm e}^{-at},
\label{xdemtsol}
\end{equation}
where $y_{\rm s}(0)$ is an arbitrary real constant which specifies any
particular solution of the pertinent homogeneous DE.

By the TLT method one obtains the solution
\begin{eqnarray}
y(t)&=&L^{-1}\Biggl\lbrace
{\hat x\omega\over s^2+\omega^2}\cdot{1\over s+a}\Biggr\rbrace
+L^{-1}\Biggl\lbrace{y(0)\over s+a}\Biggr\rbrace\nonumber\\
&=&{\hat x\over a^2+\omega^2}(\omega{\rm e}^{-at}+a\sin\omega t
-\omega\cos\omega t)+y(0){\rm e}^{-at};\hbox{\ \ for\ }t>0.
\label{xdetltsol}
\end{eqnarray}

The result (\ref{xdetltsol}) differs from the correct solution
(\ref{xdemtsol}) in two significant respects:

\noindent
a) As compared to the correct {\em evoked} response
\begin{equation}
y_{\rm e}(t)= {\hat x\over a^2+\omega^2}(a\sin\omega t-\omega\cos\omega t)
\label{xdemtevok}
\end{equation}
Eq.~(\ref{xdetltsol}) includes an extra decaying exponential.

\noindent
b) As compared to the correct {\em spontaneous} response
\begin{equation}
y_{\rm s}(t)=y_{\rm s}(0){\rm e}^{-at}
\label{xdemtspon}
\end{equation}
Eq.~(\ref{xdetltsol}) includes the initial value $y(0)$ of the {\em total}
response $y(t)=y_{\rm e}(t)+y_{\rm s}(t)$
instead of that of the spontaneous response alone, i.e., $y_{\rm s}(0)$.

The solutions (\ref{xdemtsol}) and (\ref{xdetltsol}) both are consistent with
the DE (\ref{xde}). Yet only (\ref{xdemtsol}) represents the DE's
general solution. This becomes apparent when one attempts formal
reconciliation of (\ref{xdetltsol}) with (\ref{xdemtsol}). 
Reconciliation can only be obtained by suitable choice of a particular
``initial state''. For instance, for $y_{\rm s}(0)=0$
Eqs.~(\ref{xdetltsol}, \ref{xdemtsol}) can only be reconciled by setting
\begin{equation}
y(0)={-\hat x\omega\over a^2+\omega^2}.
\label{ynullset}
\end{equation}
Hence, $y(0)$ is not really a free parameter
of the spontaneous response. While the consistency of (\ref{xdetltsol}) with
(\ref{xde}) suggests that the TLT-based solution is correct,
comparison with (\ref{xdemtsol}) reveals that (\ref{xdetltsol}) actually
is only a {\em particular} solution of (\ref{xde}) and that $y(0)$ is not
a freely assignable constant.

In point of fact, the correct solution of (\ref{xde})
can by L-transformation {\em not be obtained at all,} because the L-transform
of $\hat x\sin\omega t$ is just an alias for
$L\{u(t)\cdot\hat x\sin\omega t\}$; cf.\ Sect.~\ref{SSunilat}.
The evoked part of the
TLT-based solution (\ref{xdetltsol}) inevitably depicts for $t>0$ the
response to $u(t)\cdot\hat x\sin\omega t$ instead of to $\hat x\sin\omega t$.
The above observations reveal that LT actually is sensitive to the
difference between $x(t)$ and $u(t)x(t)$, i.e., if there is
$x(t)\ne u(t)x(t)$.

From a more general point of view,
it should be appreciated that any ``causal'' physical system which is governed
by a linear DE of finite order has a ``memory'' which becomes manifest in
the length of the system's impulse response. This notion suffices
to explain that the evoked response $y_{\rm e}(t)$
to an excitation function $x(t)$ in general depends on the behavior of $x(t)$
for $-\infty<t$ \cite{Terhardt1985,Terhardt1987}.
When the response is determined by LT, this fact is automatically
accounted for, namely. by the implicit bilaterality of (\ref{ult});
cf.\ Sects.~\ref{SSltbyft}, \ref{SSunilat}, \ref{SSinter}. The bilaterality
of $\varphi(t)$ is indispensable for the evoked response to be correct, i.e.,
whether or not $\varphi(t)$ is causal. L-transformation by (\ref{ult})
just provides for $\varphi(t)$ being always causal. As a consequence,
the evoked response complies with $x(t)$ only if 
$x(t)$  was \textit{a priori} defined to be causal.

\subsection{Analysis interval versus integration interval}
\label{SSinter}In Sect.~\ref{SSltbyft} it is pointed out that the t-domain interval
which actually is encompassed both by Fourier- and Laplace-transformation 
invariably is infinite, i.e., $-\infty<t<+\infty$. This interval is termed
the respective transformation's \textit{analysis interval}; it has to be
distinguished from the 
\textit{integration interval}, i.e., the interval of $t$ which is determined
by the transformation-integral's limits, e.g. in Eqs.~(\ref{ult}) and
(\ref{ultfour}). The relationship between these two kinds of interval
can be made apparent by regressing to \textit{discrete} Fourier- and
Laplace-transformation. The discrete transformations can be deduced from
the Fourier-series representation of a real function $f(t)$, i.e.,
\begin{equation}
f(t)=a_0+\sum_{n=1}^\infty a_{\rm n}\cos\omega_{\rm n}t
+\sum_{n=1}^\infty b_{\rm n}\sin\omega_{\rm n}t;\hbox{\ \ }n\hbox{\ integer;}
\label{fourser}
\end{equation}
where
\begin{equation}
a_0={1\over T}\int_{-T/2}^{+T/2}f(t){\ \rm d}t;
\label{fourserc}
\end{equation}
\begin{equation}
a_{\rm n}={2\over T}\int_{-T/2}^{+T/2}f(t)\cos\omega_{\rm n}t{\ \rm d}t;
\label{fourseran}
\end{equation}
\begin{equation}
b_{\rm n}={2\over T}\int_{-T/2}^{+T/2}f(t)\sin\omega_{\rm n}t{\ \rm d}t;
\label{fourserbn}
\end{equation}
\begin{equation}
\omega_{\rm n}={2\pi n\over T}.
\label{omegan}
\end{equation}
The function $f(t)$ such represented is periodic with the
period length $T=2\pi/\omega_1$.

This set of formulas can be regarded and employed as a complementary pair of
\textit{transformations} \cite{Terhardt1985}.
By Eqs.~(\ref{fourserc}-\ref{fourserbn}) a
particular section of $f(t)$ becomes transformed into an infinite set of
real coefficients \{$a_{\rm n}, b_{\rm n}$\} each of which pertains to a
discrete frequency $\omega_{\rm n}$. 
The section of $f(t)$ that extends from $-T/2$ to $+T/2$ plays the role of
an \textit{analysis interval}, $T$. 
The function $f(t)$ needs to be defined for $-T/2\le t\le+T/2$.
Inverse transformation is achieved by (\ref{fourser}). In general,
the inverse transform matches the original function $f(t)$ only for
$-T/2<t<+T/2$. By the notation chosen in
Eqs.~(\ref{fourserc}-\ref{fourserbn}) the center of the
analysis interval becomes implicitly denoted $t=0$.

To elucidate the implications of this approach it is helpful to express
the above formulas in complex notation. Using Euler's formula, one obtains
from Eqs.~(\ref{fourseran}, \ref{fourserbn}) 
\begin{equation}
a_{\rm n}\pm{\rm i}b_{\rm n}
={2\over T}\int_{-T/2}^{+T/2}f(t){\rm e}^{\pm{\rm i}\omega_{\rm n}t}{\ \rm d}t;
\hbox{\ \ }n=1,2,\ldots
\label{cplxcoef}
\end{equation}
By definition of the discrete Fourier transform
\begin{equation}
F(\omega_{\rm n},T)={T\over 2}(a_{\rm n}-{\rm i}b_{\rm n})
\label{dfours}
\end{equation}
one obtains from (\ref{cplxcoef})
\begin{equation}
F(\omega_{\rm n},T)
=\int_{-T/2}^{+T/2}f(t){\rm e}^{-{\rm i}\omega_{\rm n}t}{\ \rm d}t.
\label{dfourt}
\end{equation}
Equation (\ref{fourser}) can be expressed in the form
\begin{equation}
f(t)=a_0+{\textstyle{1\over2}}\sum_{n=1}^\infty(a_{\rm n}
-{\rm i}b_{\rm n}){\rm e}^{{\rm i}\omega_{\rm n}t}
+{\textstyle{1\over2}}\sum_{n=1}^\infty(a_{\rm n}
+{\rm i}b_{\rm n}){\rm e}^{-{\rm i}\omega_{\rm n}t}.
\label{fourscplx}
\end{equation}
Observing that by (\ref{fourserc}, \ref{dfourt}) there holds $a_0=F(0,T)/T$,
one obtains from (\ref{fourscplx}) and (\ref{dfours})
for $n=\ldots,-1,0,+1,\ldots$
the formula for inverse transformation
\begin{equation}
\varphi(t)={1\over T}\sum_{n=-\infty}^{+\infty}
F(\omega_{\rm n},T)\cdot{\rm e}^{{\rm i}\omega_{\rm n}t}
=f(t)\hbox{\ \ for\ }-T/2<t<+T/2.
\label{difourt}
\end{equation}

The transformation (\ref{dfourt}) is governed by, and confined to, the
analysis interval $T$. This interval can either be regarded as preset or
as being determined by the choice of $\omega_1$ and thus of the
spacing of the discrete frequencies $\omega_{\rm n}$, cf.\ (\ref{omegan}).
As the center of the analysis interval is for convenience denoted $t=0$
there is only one way
in which the limits of the integral (\ref{dfourt}) may become different from
$\pm T/2$, namely, by $f(t)$ being null within a sub-interval of $T$ that
borders $t=-T/2$ and/or $t=+T/2$.
\textit{Vice versa}, when one or both of the integral's limits are
different from $\pm T/2$ this invariably indicates that the function $f(t)$
contains a null interval.
In particular, when $f(t)$ is replaced with the causal function
$f_{\rm c}(t)=u(t)f(t)$ the discrete Fourier transform gets expressed by
\begin{equation}
F(\omega_{\rm n},T)
=\int_{-T/2}^{+T/2}u(t)f(t){\rm e}^{-{\rm i}\omega_{\rm n}t}{\ \rm d}t
=\int_0^{T/2}f(t){\rm e}^{-{\rm i}\omega_{\rm n}t}{\ \rm d}t
=\int_0^{T/2}f_{\rm c}(t){\rm e}^{-{\rm i}\omega_{\rm nu}t}{\ \rm d}t.
\label{dfourtc}
\end{equation}
The difference from $-T/2$ of the integral's low limit neither affects
the analysis interval $T$ nor does it make
the transformation unilateral; it just indicates that $f(t)$ is causal.

The discrete variant of \textit{Laplace} transformation emerges from
(\ref{dfourt}) by including in the integral the factor
$\exp(-\sigma t)$ and by setting the integral's low limit to $t=0$.
This yields
\begin{equation}
L_{\rm T}\{f(t)\}=\int_0^{T/2}f(t){\rm e}^{-s_{\rm n}t}{\ \rm d}t,
\label{dlt}
\end{equation}
where $L_{\rm T}$ denotes discrete L-transformation, and
\begin{equation}
s_{\rm n}=\sigma+{\rm i}\omega_{\rm n}.
\label{dsigma}
\end{equation}
Regarding the relationship between analysis interval and integration
interval, the same applies as was just pointed out above:
The integral's (\ref{dlt}) low limit $t=0$ indicates that the function
actually transformed is causal while the low limit marks the center of the
analysis interval. 

The \textit{inverse} discrete L-transform $\varphi(t)$
is determined by the expression
\begin{equation}
\varphi(t)=L_{\rm T}^{-1}\{L_{\rm T}\{f(t)\}\}
={1\over T}\sum_{n=-\infty}^{+\infty}
L_{\rm T}\{f(t)\}\cdot{\rm e}^{s_{\rm n}t}; \hbox{\ \ }-T/2<t<+T/2.
\label{dilt}
\end{equation}

When (\ref{dlt}) is employed for transformation of a non-causal function
$f(t)$, the inverse transform $\varphi(t)$ still is defined for
$-T/2<t<+T/2$; however, validity of $\varphi(t)=f(t)$ is confined to the
interval
$0<t<T/2$. In this sense, the transformation is for non-causal functions
inconsistent. For causal functions $f(t)=f_{\rm c}(t)$ the transformation
is consistent, as there holds $\varphi(t)=f_{\rm c}(t)$ for $-T/2<t<+T/2$.

Laplace transformation as defined by (\ref{ult}, \ref{ilt}) emerges from
(\ref{dlt}, \ref{dilt}) by letting $T\to\infty$, implying that
the spacing of analysis frequencies becomes infinitesimally
small such that the sum (\ref{dilt}) becomes an integral, i.e., (\ref{ilt}).
The transformation
integral's (\ref{ult}) low limit, i.e., $t=0$, marks
the center of the infinite analysis interval. 

\subsection{Redundancy of multiplication by u(t)}
\label{SSusred}When the unit step function is defined without inclusion of the
connect function $u_0(t)$, i.e., $u(t)=0$ for $t<0$; $u(t)=1$ for
$t\ge0$, there evidently holds $u^n(t)=u(t)$ ($n=1,2,\ldots$).
Hence, multiplication of $u(t)$ by itself is redundant. As $u(t)$ 
actually includes $u_0(t)$ it must be verified that this kind
of redundancy also holds when $u(t)$ is defined according to (\ref{ustep}).

For $u^n(t)=u(t)$ to hold it is necessary and sufficient
that $u_0^n(t)=u_0(t)$.
As was outlined in Sect.~\ref{SSustep}, the unit connect function $u_0(t)$
is sufficiently characterized by saying that it is
an infinite set of real numbers \{$0\ldots 1$\} that exists at $t=0$ and 
fills the gap which otherwise exists between $u(-0)=0$ and
$u(+0)=1$. From this definition it follows that any integer power of
$u_0(t)$ is characterized by precisely the same criterion:
$u_0^n(t)$ is for any $n=2,3,\ldots$ also an infinite set of real numbers
\{0\ldots 1\}. Thus, any power of $u_0(t)$ can be renamed $u_0(t)$.
Therefore, there actually holds the identity
\begin{equation}
u^n(t)=u(t);\hbox{\ \ }n=1,2,\ldots
\label{usteppow}
\end{equation}

Another important question is whether multiplication by $u(t)$ of a
\textit{derivative} of $u(t)$, i.e., of an impulse function, is also redundant
such that there holds
\begin{equation}
u(t)u^{(n)}(t)=u(t)\delta^{(n-1)}(t)=u^{(n)}(t)=\delta^{(n-1)}(t);
\hbox{\ \ }n=1,2,\ldots
\label{udelta}
\end{equation}

The identity (\ref{udelta}) turns out to be an inevitable consequence of
the fact that the impulse functions are LT-consistent (Sect.~\ref{SSimpulse},
Eq.~(\ref{impinv})), in combination with the causality theorem (\ref{varphi}).
When in (\ref{varphi}) one lets $f(t)=\delta^{(n)}(t)$ one obtains
from (\ref{impinv}) and (\ref{varphi})
\begin{equation}
L^{-1}\{L\{\delta^{(n)}(t)\}\}=\delta^{(n)}(t)=u(t)\delta^{(n)}(t).
\label{iudelta}
\end{equation}

Thus, the identity (\ref{udelta}) actually is implied in LT's
fundamental definitions and features.
The identity can be made plausible by
characterizing $\delta^{(n)}(t)$ as a set of pseudo-functions the
type of which depends on $n$ and
which only at $t=0$ are different from 0. Then it may be concluded
that multiplication
of $\delta^{(n)}(t)$ by $u_0(t)$ does not alter the type of function
denoted $\delta^{(n)}(t)$ such that, indeed, $u(t)\delta^{(n)}(t)$ can be
renamed $\delta^{(n)}(t)$. If more mathematical rigour is desired, the
theory of distributions may be invoked.

In Sect.~\ref{SSunilat} the  
identities (\ref{usteppow}) and (\ref{udelta}) are subsumed
by the \textit{unit-step redundancy theorem} (\ref{ustepequiv}).

\subsection{The dud-theorem}
\label{SSudtheo}The dud-theorem depicts the $n$-th derivative of the ud-function
$f_{\rm ud}(t)$ such as defined by (\ref{fud}), as follows.

With regard to the definition (\ref{ustep}) of the
unit step function the first derivative can be expressed by
\begin{eqnarray}
f_{\rm ud}'(t)=[u(t)f_{\rm d}(t)]'&=&0\hbox{\ \ for\ }t<0,\nonumber\\
&=&f_{\rm d}(0)u_0'(t)\hbox{\ \ for\ }t=0,\nonumber\\
&=&f_{\rm d}'(t)\hbox{\ \ for\ }t>0.
\label{fcfirst}
\end{eqnarray}
To obtain the derivatives of higher order, it is helpful that
(\ref{fcfirst}) can be converted into a more convenient form, as follows.
Utilizing the definition (\ref{ustep}) of $u(t)$,
one can replace the third line on the right side
of (\ref{fcfirst}) with the expression
\begin{equation}
u(t)f_{\rm d}'(t)-f_{\rm d}'(0)u_0(t)\hbox{\ \ for\ }-\infty<t. 
\label{fcfirstbottom}
\end{equation}
This converts (\ref{fcfirst}) into a superposition of causal
functions which are
known to be LT-consistent. As these terms are LT-consistent they can
be replaced with the pertinent inverse L-transforms without affecting the
validity of $f_{\rm ud}'(t)$. This has the effect that the term
$f_{\rm d}'(0)u_0(t)$ gets eliminated from (\ref{fcfirstbottom}).
As a result, $f_{\rm ud}'(t)$ can be expressed in the form
\begin{equation}
f_{\rm ud}'(t)=u(t)f_{\rm d}'(t)+f_{\rm d}(0)u_0'(t)\hbox{\ \ for\ }-\infty<t.
\label{fcfirstu}
\end{equation}

By comparison of (\ref{fcfirstu}) to (\ref{fud}) the scheme
becomes apparent according to which
the second derivative emerges from the first, the third from the
second, and so on. For instance, when (\ref{fcfirstu}) is derived to
obtain the second derivative of $f_{\rm ud}(t)$, the term
$u(t)f_{\rm d}'(t)$ gets converted into
$u(t)f_{\rm d}''(t)+f_{\rm d}'(0)u_0'(t)$;
and the second term $f_{\rm d}(0)u_0'(t)$ gets converted into
$f_{\rm d}(0)u_0''(t)$. By this scheme one eventually obtains for
the $n$-th derivative the expression
\begin{equation}
f_{\rm ud}^{(n)}(t)=u(t)f_{\rm d}^{(n)}(t)+f_{\rm d}^{(n-1)}(0)u_0'(t)
+f_{\rm d}^{(n-2)}(0)u_0''(t)+\ldots+f_{\rm d}(0)u_0^{(n)}(t);
\hbox{\ }n=1,2,\ldots
\label{fudn}
\end{equation}
Equation (\ref{fudn}) is by (\ref{udtheo}) expressed as the dud-theorem.

\subsection{Derivatives of order 1/2}
\label{SSxdhalf}From the derivatives of non-integer order of causal functions
those of order $1/2$ are of particular interest \cite{Sokolov}.
Below, a number of derivatives are listed which were
determined by the gdi-theorem (\ref{gditheo}), using the table
of LT-correspondences Sect.\ \ref{SScorr}.
\begin{equation}
D^{1/2}\{u(t)\}=L^{-1}\Bigl\lbrace\sqrt{s}\cdot{1\over s}\Bigr\rbrace
=u(t)\cdot{1\over\sqrt{\pi t}}
\label{dhalfustep}
\end{equation}
\begin{equation}
D^{1/2}\Bigl\lbrace{u(t)\over\sqrt{t}}\Bigr\rbrace
=L^{-1}\Bigl\lbrace\sqrt{s}\cdot\sqrt{\pi/s}\Bigr\rbrace
=\delta(t)\cdot\sqrt{\pi}
\label{dhalfovsqrt}
\end{equation}
\begin{equation}
D^{1/2}\{u(t)\cdot\sqrt{t}\}
=L^{-1}\Bigl\lbrace\sqrt{s}\cdot{\sqrt{\pi}\over2s\sqrt{s}}\Bigr\rbrace
=u(t)\cdot{\sqrt{\pi}\over2}
\label{dhalfsqrt}
\end{equation}
\begin{equation}
D^{1/2}\{u(t)\cdot{\rm e}^{t}\}
=L^{-1}\Bigl\lbrace\sqrt{s}\cdot{1\over s-1}\Bigr\rbrace
=u(t)\cdot[1/\sqrt{\pi t}+{\rm e}^t\cdot{\rm erf}(\sqrt{t})]
\label{dhalfexp}
\end{equation}
\begin{equation}
D^{1/2}\{u(t)\cdot\ln t\}
=L^{-1}\Bigl\lbrace\sqrt{s}\cdot{-\ln s-C_{\rm E}\over s}\Bigr\rbrace
=u(t)\cdot{\ln(4t)\over\sqrt{\pi t}}
\label{dhalfln}
\end{equation}
\begin{equation}
D^{1/2}\Bigl\lbrace{u(t)\over\sqrt{t}}\cdot{\rm e}^{-a^2/(4t)}\Bigr\rbrace
=L^{-1}\{\sqrt{\pi}\cdot{\rm e}^{-a\sqrt{s}}\}
=u(t)\cdot{a\over2}\cdot t^{-3/2}\cdot{\rm e}^{-a^2/(4t)}
\label{dhalfexpat}
\end{equation} 
\begin{equation}
D^{1/2}\{u(t)\cdot J_0(2\sqrt{at})\}
=L^{-1}\Bigl\lbrace{1\over\sqrt{s}}\cdot{\rm e}^{-a/s}\Bigr\rbrace
=u(t)\cdot{\cos2\sqrt{at}\over\sqrt{\pi t}}
\label{dhalfbess}
\end{equation}

\subsection{Linear DE that includes a derivative\\
of the excitation function}
\label{SSxcirc}When the linear first-order DE (\ref{xde}) is modified into the form
\begin{equation}
ay(t)+y'(t)=x'(t),
\label{xdec}
\end{equation}
the response $y(t)$ denotes the electrical current through the 
capacitor of the electrical circuit shown in Fig.~\FigCircuit, while
$x(t)$ denotes the current exerted on the circuit by the source Q.
The constant $a$ equals $a=1/(RC)$.
\begin{figure}[h]
\begin{minipage}[t]{160mm}
\begin{minipage}[t]{24mm}
\hbox{\diagram\char002}
\end{minipage}
\hfill
\parbox[b]{131mm}{\small\baselineskip10pt Fig.\ \FigCircuit.\
Electrical circuit which is accounted for by Eq.\
(\ref{xdec}).
$x(t)$: electrical current from source Q; $\eta(t)$: voltage at
condenser C; $y(t)$: current through C
\baselineskip12pt}
\end{minipage}
\end{figure}

The DE (\ref{xdec}) is of the type which by Doetsch was banished from his
theory because it includes a derivative of the excitation function.
In terms of (\ref{diffeq}) this DE corresponds to $N=1$, $M=1$.

The new method, i.e., the formula (\ref{desolg}), provides for an algorithmic,
straightforward 
solution. Assuming, as an example, $x(t)=\hat xu(t)$, and taking into
account the parameters $N=M=1$; $a_0=a$; $a_1=1$; $b_0=0$; $b_1=1$, one obtains
\begin{eqnarray}
y(t)&=&L^{-1}\Bigl\lbrace{\hat x\over s+a}\Bigr\rbrace
+L^{-1}_{\rm d}\Bigl\lbrace{y_{\rm s}(0)\over s+a}\Bigr\rbrace\nonumber\\
&=&[\hat xu(t)+y_{\rm s}(0)]\cdot{\rm e}^{-at}\hbox{\ \ for\ }-\infty<t.
\label{xamplloes}
\end{eqnarray}
Any particular initial state can be freely accounted for by setting
$y_{\rm s}(0)$ accordingly. The two components of the total solution are 
illustrated in Fig.~\FigLoes.
\begin{figure}[h]
\begin{minipage}[t]{160mm}
\begin{minipage}[t]{74mm}
\hbox{\diagram\char003}
\end{minipage}
\hfill
\parbox[b]{81mm}{\small\baselineskip10pt Fig.\ \FigLoes.\
The two components of the solution of
(\ref{xdec}).
Left: evoked response for $x(t)=\hat xu(t)$.
Right: spontaneous reponse. This
solution is depicted by Eq.\ (\ref{xamplloes})
\baselineskip12pt}
\end{minipage}
\end{figure}


\subsection{LT-correspondences} 
\label{SScorr}Below, LT-correspondences are listed that are used in the present article.
The $t$-domain functions are identical to those included in
customary tables of TLT, except for the factor $u(t)$, which is required
to make the correspondences valid in both directions.
Notice that the first correspondence (\ref{corrone}) holds only
in one direction.
\begin{equation}
1\Rightarrow 1/s
\label{corrone}
\end{equation}
\begin{equation}
u(t)\Leftrightarrow 1/s
\label{corrustep}
\end{equation}
\begin{equation}
\delta^{(n)}(t)=u^{(n+1)}(t)=u_0^{(n+1)}(t)
\Leftrightarrow s^n\hbox{\ \ \ \ }(n=0,1,\ldots)
\label{corrdelta}
\end{equation}
\begin{equation}
u(t)\cdot t^n\Leftrightarrow{n!\over s^{n+1}}\hbox{\ \ \ }(n=0,1,\ldots)
\label{corrtn}
\end{equation}
\begin{equation}
u(t)\cdot{\rm e}^{-at}\Leftrightarrow{1\over s+a}
\label{correxp}
\end{equation}
\begin{equation}
{u(t)\over a}(1-{\rm e}^{-at})\Leftrightarrow{1\over s(s+a)}
\label{correxpi}
\end{equation}
\begin{equation}
u(t)\cdot\sin\omega t\Leftrightarrow{\omega\over s^2+\omega^2}
\label{corrsin}
\end{equation}
\begin{equation}
u(t)\cdot\cos\omega t\Leftrightarrow{s\over s^2+\omega^2}
\label{corrcos}
\end{equation}
\begin{equation}
{u(t)\over a^2+\omega^2}(\omega{\rm e}^{-at}+a\sin\omega t-\omega\cos\omega t)
\Leftrightarrow{\omega\over(s+a)(s^2+\omega^2)}
\label{corrsinsol}
\end{equation}
\begin{equation}
{u(t)\over\sqrt{t}}\Leftrightarrow\sqrt{\pi\over s}
\label{corrbst}
\end{equation}
\begin{equation}
u(t)\cdot\sqrt{t}\Leftrightarrow{\sqrt{\pi/s}\over2s}
\label{corrst}
\end{equation}
\begin{equation}
u(t)\cdot{\rm e}^t\cdot{\rm erf}(\sqrt{t})\Leftrightarrow{1\over(s-1)\sqrt{s}}
\label{correerf}
\end{equation}
\begin{equation}
u(t)\cdot\ln t\Leftrightarrow{-\ln s-C_{\rm E}\over s}
\hbox{\ \ \ \ }(C_{\rm E}=0.577215\ldots)
\label{corrln}
\end{equation}
\begin{equation}
u(t)\cdot{\ln t\over\sqrt{t}}\Leftrightarrow
-\sqrt{\pi\over s}\cdot(\ln 4s+C_{\rm E})
\label{corrlnovsrt}
\end{equation}
\begin{equation}
u(t)\cdot t^r\Leftrightarrow{\Gamma(r+1)\over s^{r+1}}
\hbox{\ \ \ \ }(r\in\Re;\hbox{\ } r>-1)
\label{corrtr}
\end{equation}
\begin{equation}
{u(t)\over\sqrt{\pi t}}\cdot{\rm e}^{-a^2/(4t)}
\Leftrightarrow{1\over\sqrt{s}}\cdot{\rm e}^{-a\sqrt{s}}
\label{correxpata}
\end{equation}
\begin{equation}
u(t)\cdot a\cdot t^{-3/2}\cdot{\rm e}^{-a^2/(4t)}\Leftrightarrow
2\sqrt{\pi}\cdot{\rm e}^{-a\sqrt{s}}
\label{correxpatb}
\end{equation}
\begin{equation}
{u(t)\over\sqrt{\pi t}}\cdot\cos2\sqrt{at}\Leftrightarrow
{1\over\sqrt{s}}\cdot{\rm e}^{-a/s}
\label{corrcossqrt}
\end{equation}
\begin{equation}
u(t)\cdot J_0(2\sqrt{at})\Leftrightarrow{1\over s}\cdot{\rm e}^{-a/s}
\label{corrbess}
\end{equation}

\end{appendix}



\vfill
\baselineskip8pt
\begin{thebibliography}{99}
\begin{small}
\parskip0pt
\bibitem{Doetsch1937}
Doetsch, G.: {\it Theorie und Anwendung der Laplace-Transformation.}
Julius Springer, Berlin 1937
\bibitem{Doetsch1950}
Doetsch, G.: {\it Handbuch der La\-place-Trans\-for\-ma\-tion. I.
The\-o\-rie der La\-place-Trans\-for\-ma\-tion.} 1st ed.,
Birkh\"auser, Basel 1950
\bibitem{Doetsch1955}
Doetsch, G.: {\it Handbuch der Laplace-Transformation. II. Anwendungen der
Laplace-Transformation.} 1st ed., Birkh\"auser, Basel 1955
\bibitem{Doetsch1970}
Doetsch, G.: {\it Einf\"uhrung in Theorie und Anwendung der
Laplace-Transformation.} 2nd ed.,
Birkh\"auser, Basel Stuttgart 1970
\bibitem{Doetsch1974}
Doetsch, G.:
{\it Introduction to the Theory and Application of the Laplace Transformation.}
(transl.\ by Walter Nader).
Springer, Berlin Heidelberg New York 1974
\bibitem{Connor2004}
O'Connor, J.J., Robertson, E.F.: \textit{Gustav Doetsch}.\\
http://www.gap-system.org/\ $\tilde{}\ $history/Mathematicians/Doetsch.html
(Nov.\ 2004)
\bibitem{Oppenheim}
Oppenheim, A.V., Willsky, A.S.: {\it Signals and Systems.} 2nd ed.,
Prentice Hall, Upper Saddle River, NJ 1997
\bibitem{PhillipsParr}
Phillips, C.L., Parr, J.M., Riskin, E.A.:
{\it Signals, Systems, and Transforms.}
3rd ed., Prentice-Hall, Upper Saddle River, NJ 2002
\bibitem{Remmert1999}
Remmert, V.R.: Mathematicians at war; power struggles in Nazi
Germany's mathematical community: Gustav Doetsch and Wilhelm S\"uss.
Revue d'histoire des math\'ematiques {\bf 5} (1999), 7-59  
\bibitem{Sokolov}
Sokolov, I.M., Klafter, J., Blumen, A.: Fractional kinetics. Physics Today
{\bf 55} (Nov.~2002), 48-54
\bibitem{Terhardt1985}
Terhardt, E.: Fourier-transformation of time signals: Conceptual revision.
Acustica {\bf 57} (1985), 242-256
\bibitem{Terhardt1986}
Terhardt, E.: Ableitungsregel der Laplace-Transformation und
Anfangswertproblem. ntz-Archiv {\bf 8} (1986), 39-43
\bibitem{Terhardt1987}
Terhardt, E.: Evaluation of linear-system responses by
Laplace-transformation: Critical review and revision of method.
Acustica {\bf 64} (1987), 61-72
\bibitem{Weisstein}
Weisstein, E.W.: Laplace Transform. In: {\it MathWorld -- A Wolfram
Web Resource.}
http://mathworld.wolfram.com/LaplaceTransform.html (July 2005)
\end{small}
\end{thebibliography}
\end{document}